\documentclass[leqno,a4paper,12pt]{article}
\usepackage{amsmath,amsthm}

\newtheorem{Thm}{\indent Theorem}[section]
\newtheorem{Prop}[Thm]{\indent Proposition}
\newtheorem{Lem}[Thm]{\indent Lemma}
\newtheorem{Cor}[Thm]{\indent Corollary}
\newtheorem{Var}[Thm]{\indent Variant}
\newtheorem{Prin}[Thm]{\indent Principle}
\newtheorem{Comp}[Thm]{\indent Complement}

\theoremstyle{definition}
\newtheorem{Def}[Thm]{\indent Definition}
\newtheorem{Rem}[Thm]{\indent Remark}
\newtheorem{Ex}[Thm]{\indent Example}

\newtheorem{Constr}[Thm]{\indent Construction}
\newtheorem{Prob}[Thm]{\indent Problem}
\def\qed{{\hskip0pt\unskip\unskip\nobreak\hfil\penalty50
          \hskip1em\hbox{}\nobreak\hfil
          {\bf q.e.d.}%
          \parfillskip=0pt\finalhyphendemerits=0
          \par}\medskip}

\newenvironment{Proof}
               {{\it Proof.}\quad}
               {\qed}

\newenvironment{Proofof}[1]
               {{\it Proof of #1.}\quad}
               {\qed}

\newcommand{\Prime}{\kern3\fontdimen1\font$'$\kern-7\fontdimen1\font}

\long\def\forget#1{}

\long\def\beginSIDEREMARK#1\endSIDEREMARK
    {{\par\bigskip\advance\leftskip by 2cm
                  \advance\rightskip by -2cm\noindent
      {\bf Our own side remark:} #1
      \par\bigskip\noindent}}
\long\def\beginFORGET#1\endFORGET{#1}
\long\def\beginFORGET#1\endFORGET{}

\def\?{\ ???\ \immediate\write16{}%
\immediate\write16{Warning: There was still a question mark . . . }%
\immediate\write16{}}

\usepackage{amsmath}
\usepackage{exscale}
\usepackage{amssymb}

\newcommand{\BA}{{\mathbb{A}}}

\newcommand{\BC}{{\mathbb{C}}}
\newcommand{\BD}{{\mathbb{D}}}
\newcommand{\BE}{{\mathbb{E}}}

\newcommand{\BQ}{{\mathbb{Q}}}

\newcommand{\BZ}{{\mathbb{Z}}}

\newcommand{\FZ}{{\mathfrak{Z}}}

\newcommand{\CC}{{\cal C}}
\newcommand{\CD}{{\cal D}}

\newcommand{\CV}{{\cal V}}
\newcommand{\CW}{{\cal W}}

\newcommand{\ch}{\mathop{\rm CH}\nolimits}
\newcommand{\Spec}{\mathop{{\bf Spec}}\nolimits}

\newcommand{\imm}{\mathop{{\rm im}}\nolimits}

\newcommand{\End}{\mathop{\rm End}\nolimits}
\newcommand{\Ext}{\mathop{\rm Ext}\nolimits}
\newcommand{\GL}{\mathop{\rm GL}\nolimits}

\newcommand{\Hom}{\mathop{\rm Hom}\nolimits}

\newcommand{\Sym}{\mathop{\rm Sym}\nolimits}
\newcommand{\coker}{\mathop{\rm coker}\nolimits}

\newcommand{\loccit}{[loc.$\;$cit.]}

\def\halb{\frac{1}{2}}

\def\id{{\rm id}}

\newbox\mybox
\def\arrover#1{\mathrel{
       \setbox\mybox=\hbox spread 1.4em{\hfil$\scriptstyle#1$\hfil}
       \vbox{\offinterlineskip\copy\mybox
             \hbox to\wd\mybox{\rightarrowfill}}}}
\def\larrover#1{\mathrel{
       \setbox\mybox=\hbox spread 1.4em{\hfil$\scriptstyle#1$\hfil}
       \vbox{\offinterlineskip\copy\mybox
             \hbox to\wd\mybox{\leftarrowfill}}}}

\def\ontoover#1{\mathrel{
       \setbox\mybox=\hbox spread 1.4em{\hfil$\scriptstyle#1$\hfil}
       \vbox{\offinterlineskip\copy\mybox
             \hbox to\wd\mybox{\rightarrowfill\hskip-2.8mm
                               $\rightarrow$}}}}
\def\leftontoover#1{\mathrel{
       \setbox\mybox=\hbox spread 1.4em{\hfil$\scriptstyle#1$\hfil}
       \vbox{\offinterlineskip\copy\mybox
             \hbox to\wd\mybox{$\leftarrow$\hskip-2.8mm
                               \leftarrowfill}}}}
\def\longto{\longrightarrow}
\def\into{\hookrightarrow}

\def\longonto{\ontoover{\ }}
\def\isoto{\arrover{\sim}}

\def\longinto{\lhook\joinrel\longrightarrow}

\usepackage[curve,matrix,arrow,cmtip]{xy}
\NoComputerModernTips

\def\myxymessage{\def\messagetext
   {Here an xy-pic diagram was omitted to speed up compilation . . . }
   \immediate\write16{\messagetext}
   \hbox{\bf \messagetext}}
\def\filxymatrix#1{\myxymessage}
\def\filxyarray#1{\myxymessage}
\newdir^{ (}{{}*!/-3pt/\dir^{(}}
\newdir_{ (}{{}*!/-3pt/\dir_{(}}
\newdir^{ )}{{}*!/+3pt/\dir^{)}}
\newdir_{ )}{{}*!/+3pt/\dir_{)}}

\def\rscript#1{\hbox to 0pt{$\scriptstyle#1$\hss}}

\let\oldbullet\bullet
\def\bullet{{\mathchoice{\oldbullet}%
                        {\oldbullet}%
                        {\scriptscriptstyle\oldbullet}%
                        {\oldbullet}}}

\newcommand{\argdot}{{\;\bullet\;}}%Punkt als Platzhalter fuer Argumente
\newcommand{\bS}{\mathop{\overline{S}}\nolimits}
\newcommand{\bT}{\mathop{\overline{T}}\nolimits}
\newcommand{\bU}{\mathop{\overline{U}}\nolimits}
\newcommand{\bX}{\mathop{\overline{X}}\nolimits}
\newcommand{\bY}{\mathop{\overline{Y}}\nolimits}
\newcommand{\dX}{\mathop{\partial \overline{X}}\nolimits}
\newcommand{\dY}{\mathop{\partial \overline{Y}}\nolimits}
\newcommand{\uh}{\mathop{\underline{h}}\nolimits}

\newcommand{\zeq}{\mathop{z_{equi}}\nolimits}
\newcommand{\CHeffM}{\mathop{CHM^{eff}(k)}\nolimits}
\newcommand{\CHSeffM}{\mathop{CHM^{s,eff}(S)}\nolimits}
\newcommand{\CHM}{\mathop{CHM(k)}\nolimits}
\newcommand{\CHSM}{\mathop{CHM^s(S)}\nolimits}
\newcommand{\CHTM}{\mathop{CHM^s(T)}\nolimits}
\newcommand{\CHUM}{\mathop{CHM^s(U)}\nolimits}

\newcommand{\CHSQM}{\mathop{CHM^s(S)_F}\nolimits}
\newcommand{\DeffgM}{\mathop{DM^{eff}_{gm}(k)}\nolimits}
\newcommand{\DeffQgM}{\mathop{\DeffgM_F}\nolimits}
\newcommand{\DeffqgM}{\mathop{\DeffgM_\BQ}\nolimits}

\newcommand{\DgM}{\mathop{DM_{gm}(k)}\nolimits}

\newcommand{\DSgM}{\mathop{DM_{gm}(S)}\nolimits}

\newcommand{\Mcgm}{\mathop{M^c}\nolimits}

\newcommand{\dMgm}{\mathop{\partial M}\nolimits}
\newcommand{\RC}{\mathop{{\bf R} C}\nolimits}
\newcommand{\Xp}{\mathop{\widetilde{X}} \nolimits}
\newcommand{\Yp}{\mathop{\widetilde{Y}} \nolimits}

\begin{document}

%%%%%%%%%%%%%%%%%%%%%%%%%%%%%%%%%%%%%%%%%%%%%%%%%%%%%%%%%%%%%%%%%%%%%%%
%
%  formatting

\hfuzz=3pt
\overfullrule=10pt                   % erzeugt schwarze Fehlerbalken

% The displayskip values were changed because \LaTeX does not react
% correctly to a \leqno: it should then use big skips, but doesn't.

\setlength{\abovedisplayskip}{6.0pt plus 3.0pt}
                               % preset 10.0pt plus 2.0pt minus 5.0pt
\setlength{\belowdisplayskip}{6.0pt plus 3.0pt}
                               % preset 10.0pt plus 2.0pt minus 5.0pt
\setlength{\abovedisplayshortskip}{6.0pt plus 3.0pt}
                               % preset 0.0pt plus 3.0pt
\setlength{\belowdisplayshortskip}{6.0pt plus 3.0pt}
                               % preset 6.0pt plus 3.0pt minus 3.0pt

\setlength{\baselineskip}{13.0pt}
                               % preset 12.0pt
\setlength{\lineskip}{0.0pt}
                               % preset 1.0pt
\setlength{\lineskiplimit}{0.0pt}
                               % preset 0.0pt

%%%%%%%%%%%%%%%%%%%%%%%%%%%%%%%%%%%%%%%%%%%%%%%%%%%%%%%%%%%%%%%%%%%%%%%
%
%  Title Page
%
%%%%%%%%%%%%%%%%%%%%%%%%%%%%%%%%%%%%%%%%%%%%%%%%%%%%%%%%%%%%%%%%%%%%%%%

\title{Boundary motive, relative motives and extensions of motives
\forget{
\footnotemark
\footnotetext{To appear in ....}
}
}
\author{\footnotesize by\\ \\
\mbox{\hskip-2cm
\begin{minipage}{6cm} \begin{center} \begin{tabular}{c}
J\"org Wildeshaus \footnote{
Partially supported by the \emph{Agence Nationale de la
Recherche}, project no.\ ANR-07-BLAN-0142 ``M\'ethodes \`a la
Voevodsky, motifs mixtes et G\'eom\'etrie d'Arakelov''. }\\[0.2cm]
\footnotesize LAGA\\[-3pt]
\footnotesize UMR~7539\\[-3pt]
\footnotesize Institut Galil\'ee\\[-3pt]
\footnotesize Universit\'e Paris 13\\[-3pt]
\footnotesize Avenue Jean-Baptiste Cl\'ement\\[-3pt]
\footnotesize F-93430 Villetaneuse\\[-3pt]
\footnotesize France\\
{\footnotesize \tt wildesh@math.univ-paris13.fr}
\end{tabular} \end{center} \end{minipage}
\hskip-2cm}
\\[2.5cm]
\forget{
{\bf Preliminary version --- not for distribution}\\[1cm]
}
}
% In the final version we might want to fix the date:
\date{May 5, 2011}
\maketitle
%\quad \\[-1.7cm]
\begin{abstract}
\noindent 
We explain the role of the boundary motive in the construction
of certain Chow motives, and of extensions of Chow motives.
Our two main examples concern proper, singular surfaces 
and fibre products of a universal elliptic curve.  \\


%\noindent
%{\bf R\'esum\'e~:} RESUME.\\
\end{abstract}

%\vfill

\bigskip
\bigskip
\bigskip

\noindent {\footnotesize Math.\ Subj.\ Class.\ (2010) numbers: 19E15
(14F25, 14F42, 14G35). 
}

\eject
\tableofcontents

\bigskip
%\vspace*{0.5cm}

%\newpage
%\include{Intro}

%%%%%%%%%%%%%%%%%%%%%%%%%%%%%%%%%%%%%%%%%%%%%%%%%%%%%%%%%%%%%%%%%%%%%%%
%
%  Introduction
%
%%%%%%%%%%%%%%%%%%%%%%%%%%%%%%%%%%%%%%%%%%%%%%%%%%%%%%%%%%%%%%%%%%%%%%%

\setcounter{section}{-1}
\section{Introduction}
\label{Intro}

%%%%%%%%%%%%%%%%%%%%%%%%%%%%%%%%

%%%%%%%%%%%%%%%%%%%%%%%%%%%%%%%%

This article contains largely extended notes
of a short series of lectures delivered during the
\emph{Ecole d'\'et\'e franco-asiatique ``Autour des motifs''},
which took place at the IHES in July~2006. 
The task which I was assigned was to explain the role 
of the \emph{boundary motive}, and I hope that
the present article will make a modest contribution to this effect. \\ 

By definition \cite{W1}, 
the boundary motive $\dMgm(X)$ of a variety $X$ over a perfect field $k$
fits into a canonical exact triangle
\[
(\ast) \quad\quad
\dMgm(X) \longto M(X) \longto \Mcgm(X) \longto \dMgm(X)[1]
\]
in the category $\DeffgM$ of effective \emph{geometrical motives}.
This triangle establishes the relation of
the boundary motive to $M (X)$ and $\Mcgm (X)$, 
the \emph{motive} of $X$ and its \emph{motive with compact support}, 
respectively \cite{V}. \\

One way to explain its interest is to start with the notion
of extensions. Indeed, most of the existing attempts 
to prove the Beilinson or
Bloch--Kato conjectures on special values of $L$-functions
necessitate
the construction of extensions of (Chow) motives,
and the explicit control of their realizations (Betti, de~Rham,
\'etale...). Often, 
the source of these extensions is \emph{localization},
which expresses the motive with compact support of a 
non-compact variety $X$ as an extension of the motive of
a compactification $X^*$ by the motive of the complement $X^* - X$.
The realizations of these extensions then correspond
to cohomology with compact support of $X$.
This approach is clearly present e.g.\ in Harder's work on 
special values \cite{H}. \\

Thus,
given two Chow motives, one may try to use localization
to construct an extension of one by the other.
Here, we base ourselves on the principle that the given
Chow motives are ``basic'', and that the extension is
``difficult'' to obtain. 
But one may also invert the logic:
given a ``mixed'' motive, try to use localization
to construct the Chow motives used to build it up;
let us refer to this problem as ``resolution of extensions''. \\

The purpose of this article is to establish that the
boundary motive plays a role both for
the construction and for the resolution of extensions
\emph{via} localization. 
In Section~\ref{1}, we start by making precise the relation
between localization and the boundary motive. In fact, the
triangle $(*)$ turns out to be obtained by ``splicing'' the localization
triangle and its dual. 
We chose to discuss this relation first in 
the Hodge theoretic realization, and in the special
case of a complement $X$ of two points in an elliptic curve
over $\BC$ (Examples~\ref{1A}, \ref{1C} and \ref{1E}),
and deduce from that discussion the general picture
in Hodge theory (Theorems~\ref{1F} and \ref{1G}), 
concerning compactifications
of a fixed variety $X$ over $\BC$. We observe in particular 
(Corollary~\ref{1H}) that when $X$ is smooth, 
then any smooth compactification induces a
\emph{weight filtration} on the boundary cohomology of $X$, i.e.,
on the Hodge realization of the boundary motive. \\

In order to formulate the motivic analogues of these results,
we need the right notion of weights for motives.
It turns out that this notion is given by \emph{weight structures},
as recently introduced and studied by Bondarko \cite{Bo}.
We review the definition, and the basic properties of weight
structures, including their application to motives (Theorem~\ref{1L}):
according to Bondarko, there is a canonical such structure on
the triangulated category $\DeffgM$, and its \emph{heart} equals
the category $\CHeffM$ of effective Chow motives.
The motivic analogue of Corollary~\ref{1H} holds: according to
Corollary~\ref{1P}, any smooth compactification of a fixed variety
$X$ which is smooth over $k$ induces a weight filtration on $\dMgm(X)$. \\

Then we try to invert this process (hoping for this inversion to allow us to
resolve extensions). The precise statement is given
in Theorem~\ref{1R}, which states that for fixed $X$,
there is a canonical bijective correspondence
(discussed at length in Construction~\ref{1Q}) between isomorphism
classes of two types of objects: (1)~weight filtrations on $\dMgm(X)$,
and (2)~certain effective Chow motives $M_0$ through which the morphism
$M(X) \to \Mcgm(X)$ factors. 
An analogous statement (Variant~\ref{1V}) holds for direct
factors of $\dMgm(X)$, $M(X)$, and $\Mcgm(X)$, provided that they
are images of an idempotent endomorphism of the whole exact triangle $(*)$.  
In this correspondence,
the passage to isomorphism classes
cannot be avoided because of the necessity to \emph{choose}
cones of certain morphisms in the triangulated category $\DeffgM$.
This causes (at least) one important problem, namely the  
lack of functoriality of the representatives of the isomorphism
classes. In order to obtain functoriality, Construction~\ref{1Q}
thus needs to be rigidified. \\

In the rest of Section~\ref{1}, we describe the approach from \cite{W3} 
to rigidification, hence functoriality. It is based on 
the notion of motives
\emph{avoiding certain weights}. If a direct factor $\dMgm(X)^e$ of $\dMgm(X)$
is without weights $-1$ and $0$, then an effective Chow
motive $M_0$ is canonically and functorially defined  
(Complement~\ref{1W}). 
Given the nature of the realizations of $M_0$, it is natural to
call it the \emph{$e$-part of the interior motive} of 
$X$. Its main properties are established in \cite[Sect.~4]{W3}.
Note however (Problem~\ref{1U}) that the above condition
on absence of weights is never satisfied for the whole of $\dMgm(X)$ ---
unless $\dMgm(X) = 0$. 
In order to make this approach work, we thus
need an idempotent endomorphism
$e$ of the exact triangle $(\ast)$, giving rise to a direct factor
\[
\dMgm(X)^e \longto M(X)^e \longto \Mcgm(X)^e \longto \dMgm(X)^e[1] \; .
\]  

Section~\ref{2b} shows how the theory of \emph{smooth relative 
Chow motives} can be employed to construct endomorphisms
of the exact triangle $(\ast)$. 
Fix a base scheme $S$, which is smooth over $k$. Theorem~\ref{2bB}
establishes the existence of a functor from the category of smooth relative
Chow motives over $S$ to the category of exact triangles in
$\DeffgM$. On objects, it is given by mapping a proper, smooth $S$-scheme
$X$ to the exact triangle 
\[
\dMgm(X) \longto M(X) \longto \Mcgm(X) \longto \dMgm(X)[1] \; .
\]
We should mention that as far as the $M(X)$-component
is concerned, the functoriality statement from Theorem~\ref{2bB}
is just a special feature of results by D\'eglise \cite{Deg2},
Cisinski--D\'eglise 
\cite{CDeg} and Levine \cite{L3}
(see Remarks~\ref{2bC} and \ref{2bGa} for details).
However, the application of the results from \loccit \ to the functor $\dMgm$
is not obvious. This is one of the reasons
why we follow an alternative
approach. It is based on a relative version
of \emph{moving cycles} \cite[Thm.~6.14]{W1}.
This also explains why we are forced to suppose the base field $k$
to admit a strict version of resolution of singula\-ri\-ties. 
Theorem~\ref{2bBa} and Corollary~\ref{2bK}
then analyze the behaviour of the functor
from Theorem~\ref{2bB} under change of the base $S$. 
Another reason for us to choose a cycle theoretic approach
was that it becomes then easier to keep track of the 
correspondences on $X \times_k X$ commuting with our constructions. 
Our main application (Example~\ref{2bL}) thus
concerns correspondences ``of Hecke type''
yielding endomorphisms of the exact triangle $(*)$. \\

In Section~\ref{2c}, we apply these principles to Abelian
schemes. More precisely, the main result of \cite{DM} on
the Chow--K\"unneth 
decomposition of the relative motive of an Abelian scheme $A$
over $S$ (recalled in Theorem~\ref{2cA})
yields canonical projectors in the 
relative Chow group. 
Given our analysis from Section~\ref{2b}, it follows that they act
idempotently on the exact triangle  
\[
\dMgm(A) \longto M(A) \longto \Mcgm(A) \longto \dMgm(A)[1] \; .
\]

In Sections~\ref{2} and \ref{3}, we discuss two examples. 
Section~\ref{2} concerns normal, proper surfaces $X^*$. We first
recall the basic construction of the \emph{intersection motive}
$M^{!*} (X^*)$ of $X^*$, following previous work of Cataldo and
Miglio\-rini \cite{CM}, and review some of the material
from \cite{W4}. In particular (Proposition~\ref{2C}), we recall
that $M^{!*} (X^*)$ is co- and contravariantly functorial
under finite morphisms of proper surfaces. We then analyze the
precise relation to the weight filtration of the boundary motive
of a dense, open sub-scheme $X \subset X^*$, which is smooth over $k$
(Theorem~\ref{2D}), following the lines of Construction~\ref{1Q}.
We finish the section with a discussion of
the case of Baily--Borel compactifications of
Hilbert surfaces. We recall, following \cite[Sect.~6 and 7]{W4},
that localization allows to construct non-trivial extensions
of a certain Artin motive by a direct factor of $M^{!*} (X^*)$.
Using Proposition~\ref{2C}, we then establish
stability of $M^{!*} (X^*)$ under the
correspondences ``of Hecke type'' constructed in Example~\ref{2bL}. \\[0.1cm]
 
In Section~\ref{3}, we discuss fibre products of the universal
elliptic curve over the modular curve of level $n \ge 3$.
We review some of the material from \cite{Sch} and \cite[Sect.~3 and 4]{W3}.
Notably (Proposition~\ref{3C}), 
we recall that in this geometrical setting, the condition
from Complement~\ref{1W} on the absence of weights $-1$ and $0$
in the boundary motive is satisfied. Thus, the interior motive
can be defined. The new ingredient is Example~\ref{3D}, where we 
use rigidity of our construction to give a proof ``avoiding compactifications''
of equivariance of the interior motive
under the correspondences ``of Hecke type''. \\
 
As mentioned above, this article is primarily intended to be a
general introduction to the construction and to the applications
of boundary motives. For many details of the proofs, we shall
refer to our earlier articles \cite{W1} and \cite{W3}. 
Let us however indicate that various parts of this paper
discuss original constructions. This is true in particular for 
Section~\ref{2b} (on relative motives and 
functoriality), including our study of Hecke equivariance. 
We expect these constructions to be of interest in other contexts
than those discussed in Sections~\ref{2} and \ref{3}. \\

For further developments of the theory of boundary motives 
and their applications to
special classes of algebraic varieties and to their associated
motives, in particular to the motives of
Shimura varieties, we refer also to \cite{W2,W5}. \\

Part of this work was done while I was enjoying a 
\emph{modulation de service pour les porteurs de projets de recherche},
granted by the \emph{Universit{\'e} Paris~13}, and during a stay at the 
\emph{Universit\"at Z\"urich}. I am grateful to both institutions.
I wish to thank the organizers of \emph{Autour des motifs}
for the invitation to Bures-sur-Yvette, and
J.~Ayoub, F.~D\'eglise, D.~H\'ebert, B.~Kahn, F.~Lecomte 
and M.~Levine for useful discussions and comments. Special thanks
go to J.-B.~Bost for insisting on this article to be written,
and for his helpful suggestions to improve an earlier version. 

\newpage

{\bf Notation and conventions}: $k$ denotes a fixed perfect
base field, $Sch/k$  
the category of separated schemes of finite type over $k$, and  
$Sm/k \subset Sch/k$
the full sub-category of objects which are smooth over $k$.
When we assume $k$ to admit resolution of singularities,
then it will be 
in the sense of \cite[Def.~3.4]{FV}:
(i)~for any $X \in Sch/k$, there exists an abstract blow-up $Y \to X$
\cite[Def.~3.1]{FV} whose source $Y$ is in $Sm/k$,
(ii)~for any $X, Y \in Sm/k$, and any abstract blow-up $q : Y \to X$,
there exists a sequence of blow-ups 
$p: X_n \to \ldots \to X_1 = X$ with smooth centers,
such that $p$ factors through $q$. 
We say that $k$ admits strict resolution of singularities,
if in (i), for any given dense open subset $U$
of the smooth locus of $X$,
the blow-up $q: Y \to X$ can be chosen to be an isomorphism above $U$,
and such that arbitrary intersections of
the irreducible components of the complement $Z$ of $U$ in $Y$  
are smooth (e.g., $Z \subset Y$ a normal
crossing divisor with smooth irreducible components). \\   

As far as motives are concerned,
the notation of this paper follows
that of \cite{V}. We refer to Levine's lecture notes (this volume) for a 
review of this notation, and in particular,
of the definition
of the categories $\DeffgM$ and $\DgM$ of (effective) geometrical
motives over $k$,
and of the motive $M(X)$ and the motive with compact support $\Mcgm(X)$
of $X \in Sch / k$. Let $F$ be a commutative flat $\BZ$-algebra,
i.e., a commutative unitary ring whose additive group is without torsion.
The notation $\DeffgM_F$ and $\DgM_F$ stands for the 
$F$-linear analogues of $\DeffgM$ and $\DgM$
defined in \cite[Sect.~16.2.4 and Sect.~17.1.3]{A}. 
Simi\-larly,
let us denote by $\CHeffM$ and $\CHM$ the categories opposite to the categories
of (effective) Chow motives, and by
$\CHeffM_F$ and $\CHM_F$ the pseudo-Abelian
completion of the category
$\CHeffM \otimes_\BZ F$ and $\CHM \otimes_\BZ F$,
respectively. Using \cite[Cor.~2]{V3} (\cite[Cor.~4.2.6]{V} if $k$ 
admits resolution of singularities), we canonically identify 
$\CHeffM_F$ and $\CHM_F$ with
a full additive sub-category of $\DeffgM_F$ and $\DgM_F$, respectively.
Note in particular that with these conventions, $\CHM_\BQ$ is 
actually \emph{opposite}
to the category denoted by the same symbol in \cite{W4}.

%%% Local Variables:
%%% mode: latex
%%% TeX-master: "head"
%%% End:

\bigskip

%\include{Sec1}

%%%%%%%%%%%%%%%%%%%%%%%%%%%%%%%%%%%%%%%%%%%%%%%%%%%%%%%%%%%%%%%%%%%%%%%
%
%  Section 1
%
%%%%%%%%%%%%%%%%%%%%%%%%%%%%%%%%%%%%%%%%%%%%%%%%%%%%%%%%%%%%%%%%%%%%%%%

\section{Motivation}
\label{1}

%%%%%%%%%%%%%%%%%%%%%%%%%%%%%%%%

%%%%%%%%%%%%%%%%%%%%%%%%%%%%%%%

Let us start by recalling the geometrical interpretation
of (cup product with) the \emph{Chern class} in a very special context. 

\begin{Ex} \label{1A}
Let $E$ be an elliptic curve over the field $\BC$ of complex numbers,
and $P \in E(\BC)$ a point unequal to zero. Put $X := E - \{ 0 , P \}$,
and consider the complementary inclusions
\[
\vcenter{\xymatrix@R-10pt{
        X \ar@{^{ (}->}[r]^-{j} &
        E \ar@{<-^{ )}}[r]^-{i} &
        \{ 0 , P \} \; .
\\}}
\]
Let us prepare the reader that
the following aspect of this situation will 
be generalized in the sequel: $E$ is a smooth compactification of $X$. 
The associated long exact \emph{localization sequence} for 
(singular) cohomology with coefficients in $\BQ$
reads as follows.
\[
\vcenter{\xymatrix@R-10pt{
        0 \ar[r] &
        H^0_c (X (\BC),\BQ) \ar[r] &
        H^0 (E (\BC),\BQ) \ar[r]  &
        H^0 (\{ 0 \},\BQ) \oplus H^0 (\{ P \},\BQ) \\
         \ar[r] &
        H^1_c (X (\BC),\BQ) \ar[r] &
        H^1 (E (\BC),\BQ) \ar[r] &
        0 \\
         \ar[r] &  
        H^2_c (X (\BC),\BQ) \ar[r] &
        H^2 (E (\BC),\BQ) \ar[r] &
        0  
\\}}
\]
It shows that $H^0_c (X (\BC),\BQ) = 0$, that
\[
H^2_c (X (\BC),\BQ) \isoto H^2 (E (\BC),\BQ) \; ,
\]
and, most interestingly, that $H^1_c (X (\BC),\BQ)$ is a Yoneda
one-extension of $H^1 (E (\BC),\BQ)$ by the cokernel of
\[
i^*: H^0 (E (\BC),\BQ) \longto 
H^0 (\{ 0 \},\BQ) \oplus H^0 (\{ P \},\BQ) \; .
\]
The Hodge structures on the three groups $H^0 (E (\BC),\BQ)$, 
$H^0 (\{ 0 \},\BQ)$ and $H^0 (\{ P \},\BQ)$ all equal 
$\BQ(0)$,
and under these identifications, $i^*$ corresponds to the diagonal
embedding 
\[
\Delta: \BQ(0) \longinto \BQ(0) \oplus \BQ(0) \; .
\]
We identify its cokernel with $\BQ(0)$ by \emph{choosing} 
the class of 
\[
(-1,1) \in H^0 (\{ 0 \},\BQ) \oplus H^0 (\{ P \},\BQ)
\]
as its generator. The one-extension
then takes the form
\[
0 \longto \BQ(0) \longto H^1_c (X (\BC),\BQ)
\longto H^1 (E (\BC),\BQ) \longto 0 \; .
\]
Let us denote it by $\BE xt_P$. Defining $\BE xt_0$ to be the trivial
extension, we thus get a map 
\[
\BE xt: E(\BC) \longto 
\Ext^1 \bigl( H^1 (E (\BC),\BQ) , \BQ(0) \bigr) \; , \;
P \longmapsto \BE xt_P 
\]
($\Ext^1 :=$ the group of one-extensions in the category
of mixed $\BQ$-Hodge structures), 
which can be checked to be a morphism of groups.
\end{Ex}

An analogous construction is possible for $\ell$-adic cohomology
and elliptic curves over a field of characteristic unequal to $\ell$.
The one-extensions then take place in the category of modules over
the absolute Galois group. 
Note that in the context considered in Example~\ref{1A},
the morphism $\BE xt$ induces an isomorphism 
\[
E(\BC) \otimes_\BZ \BQ \isoto 
\Ext^1 \bigl( H^1 (E (\BC),\BQ) , \BQ(0) \bigr) \; .
\]
Here is what we would like the reader to recall from the above.

\begin{Prin} \label{1B}
Localization potentially leads to interesting extensions of Hodge structures
or Galois modules.
\end{Prin}

Actually, Principle~\ref{1B} admits a more general version,
where we replace ``localization'' by ``the formalism of six operations''.
In the sequel of this article, we shall however
concentrate on localization and its dual. 

\begin{Ex} \label{1C}
We continue to consider the situation from Example~\ref{1A}. \\[0.1cm]
(a)~The long exact sequence dual to the localization sequence associated to 
\[
\vcenter{\xymatrix@R-10pt{
        X \ar@{^{ (}->}[r]^-{j} &
        E \ar@{<-^{ )}}[r]^-{i} &
        \{ 0 , P \}
\\}}
\]
($X = E - \{ 0 , P \}$ as before) will be referred to 
as the \emph{co-localization sequence}.
It shows that
\[
H^0 (E (\BC),\BQ) \isoto H^0 (X (\BC),\BQ) \; ,
\]
that $H^2 (X (\BC),\BQ) = 0$, 
and that $H^1 (X (\BC),\BQ)$ is a Yoneda
one-extension 
\[
0 \longto H^1 (E (\BC),\BQ) \longto H^1 (X (\BC),\BQ)
\longto \BQ(-1) \longto 0 
\]
(with the identifications dual to the one used in Example~\ref{1A}). \\[0.1cm]
(b)~Let us now compare cohomology and cohomology with compact support of $X$.
This comparison is expressed by a third long exact sequence,
which we shall refer to as the \emph{boundary sequence}.  
\[
\vcenter{\xymatrix@R-10pt{
          &
        H^0_c (X (\BC),\BQ) \ar[r] &
        H^0 (X (\BC),\BQ) \ar[r]  &
        \partial H^0 (X (\BC),\BQ) \\
         \ar[r] &
        H^1_c (X (\BC),\BQ) \ar[r] &
        H^1 (X (\BC),\BQ) \ar[r] &
        \partial H^1 (X (\BC),\BQ) \\
         \ar[r] &  
        H^2_c (X (\BC),\BQ) \ar[r] &
        H^2 (X (\BC),\BQ) &
\\}}
\]
The third term in this sequence is \emph{boundary cohomology} of $X$, defined
as cohomology of $E(\BC)$ with coefficients in the complex 
$i^*j_* \BQ_X$. Here, the symbol $j_* \BQ_X$
denotes the total direct image of $\BQ_X$ under $j$; 
following the convention used in \cite{BBD},
we drop the letter ``$R$'' from our notation.
The morphism from $H^1_c (X (\BC),\BQ)$ to 
$H^1 (X (\BC),\BQ)$ factors over $H^1 (E (\BC),\BQ)$.
Therefore, by what was said before, we see that the boundary sequence
is obtained by splicing the two sequences
\[
\vcenter{\xymatrix@R-10pt{
        0 \ar@{=}[r] &
        H^0_c (X (\BC),\BQ) \ar[r] &
        H^0 (X (\BC),\BQ) \ar[r]  &
        \partial H^0 (X (\BC),\BQ) \\
         \ar[r] &
        H^1_c (X (\BC),\BQ) \ar[r] &
        H^1 (E (\BC),\BQ) \ar[r] &
        0
\\}}
\]
and
\[
\vcenter{\xymatrix@R-10pt{
        0 \ar[r] &
        H^1 (E (\BC),\BQ) \ar[r] &
        H^1 (X (\BC),\BQ) \ar[r] &
        \partial H^1 (X (\BC),\BQ) \\
         \ar[r] &  
        H^2_c (X (\BC),\BQ) \ar[r] &
        H^2 (X (\BC),\BQ) \ar@{=}[r] &
        0    \; .
\\}}
\]
This allows us in particular to identify boundary cohomology:
\[
\partial H^0 (X (\BC),\BQ) \cong \BQ(0)^2 \quad , \quad
\partial H^1 (X (\BC),\BQ) \cong \BQ(-1)^2 \; .
\]
(c)~We claim that the above ``first half'' of the boundary sequence
\[
\vcenter{\xymatrix@R-10pt{
        0 \ar@{=}[r] &
        H^0_c (X (\BC),\BQ) \ar[r] &
        H^0 (X (\BC),\BQ) \ar[r]  &
        \partial H^0 (X (\BC),\BQ) \\
         \ar[r] &
        H^1_c (X (\BC),\BQ) \ar[r] &
        H^1 (E (\BC),\BQ) \ar[r] &
        0
\\}}
\]
equals the localization sequence
\[
\vcenter{\xymatrix@R-10pt{
        0 \ar@{=}[r] &
        H^0_c (X (\BC),\BQ) \ar[r] &
        H^0 (E (\BC),\BQ) \ar[r]  &
        H^0 (\{ 0 \},\BQ) \oplus H^0 (\{ P \},\BQ) \\
         \ar[r] &
        H^1_c (X (\BC),\BQ) \ar[r] &
        H^1 (E (\BC),\BQ) \ar[r] &
        0 
\\}}
\]
from Example~\ref{1A}. Indeed, the identification is achieved by 
the canonical isomorphism 
\[
H^0 (\{ 0 \},\BQ) \oplus H^0 (\{ P \},\BQ) 
\isoto \partial H^0 (X (\BC),\BQ) \; ,
\]
induced by the adjunction $i^* \BQ_E \to i^* j_* \BQ_X$. 
A dual statement relates the co-localization sequence from (a) 
to the ``second half''
of the boundary sequence. \\[0.1cm]
(d)~Altogether, we see that the long exact boundary sequence allows us to
recover cohomology of $E$, together with the localization and
co-localization sequences. One ``half'' of boundary cohomology
(namely $\partial H^0$) contributes to the
localization sequence, the other ``half'' (namely 
$\partial H^1$) to the co-localization sequence.
\end{Ex}

Here is what we would the like the reader to recall from the above.

\begin{Prin} \label{1D}
The boundary sequence 
allows to recover cohomology of a smooth compactification of $X$,
together with the localization and co-localization sequences.
\end{Prin}

A few precisions are necessary. First, given that $X$ is a curve,
there is only one possible choice of smooth compactification (namely $E$).
But this changes of course in higher dimensions.
Second, the ``recovery'' of the localization and co-localization
sequences from the boundary sequence seems to require a choice of 
additional data, namely a
division of boun\-dary cohomology into two ``halfs''.
In order to address both points in a satisfactory manner
(see Theorem~\ref{1F} below), we need to formalize the problem. \\

Since we wish the discussion to apply to the triangulated category of motives, 
for which no $t$-structure is available at present,
it is best placed in the context of triangulated categories.
In the context we chose to discuss, namely that of Hodge theory,
the appropriate triangulated category is the category
of \emph{algebraic $\BQ$-Hodge modules} \cite{Sa}.
We should immediately reassure readers not familiar with this theory:
for our purposes, only its formal properties (localization,
purity, proper base change,...) will be needed. Therefore,
in order to motivate what is to follow, we might 
just as well have placed ourselves
in the context of $\ell$-adic sheaves, which would allow to 
argue in a completely analogous
fashion. Readers wishing nonetheless to have a survey on 
Hodge theory at their disposal might find it useful to consult
\cite{St}.

\begin{Ex} \label{1E}
The relation to the geometric situation 
\[
\vcenter{\xymatrix@R-10pt{
        X \ar@{^{ (}->}[r]^-{j} &
        E \ar@{<-^{ )}}[r]^-{i} &
        \{ 0 , P \}
\\}}
\]
studied before is as follows.
The localization sequence from Example~\ref{1A} is the result of 
application of the cohomological functor $H^*(E(\BC),\argdot)$
to the exact \emph{localization triangle} 
\[
j_! \BQ_X(0) \longto \BQ_E(0) \longto i_* i^* \BQ_E(0)
\longto j_! \BQ_X(0)[1]
\]
of algebraic $\BQ$-Hodge modules on $E$ \cite[(4.4.1)]{Sa}.
In the same way, the co-localization sequence from Example~\ref{1C}~(a) 
is induced by the exact \emph{co-localization triangle} 
\[
i_* i^! \BQ_E(0) \longto \BQ_E(0) \longto j_* \BQ_X(0)
\longto i_* i^! \BQ_E(0)[1]
\]
of Hodge modules on $E$ \cite[(4.4.1)]{Sa}, using in addition that thanks to 
\emph{purity}, we have a canonical identification
\[
i_* i^! \BQ_E(0) \cong \BQ_E(-1)[-2] \; .
\]
Applying localization to the Hodge module $j_* \BQ_X(0)$
(or equivalently, co-localization to $j_! \BQ_X(0)$),
we obtain the exact \emph{boundary triangle}
\[
j_! \BQ_X(0) \longto j_* \BQ_X(0) \longto i_* i^* j_* \BQ_X(0)
\longto j_! \BQ_X(0)[1] \; ,
\]
which induces the boundary sequence from Example~\ref{1C}~(b).
\end{Ex}

Note that the three triangles (localization, co-localization and boundary)
exist for \emph{any} pair of complementary immersions. 
The following results from Saito's formalism
of six operations on algebraic Hodge modules \cite{Sa}.

\begin{Thm} \label{1F}
Let 
\[
\vcenter{\xymatrix@R-10pt{
        X \ar@{^{ (}->}[r]^-{j} &
        \bX \ar@{<-^{ )}}[r]^-{i} &
        D
\\}}
\]
be complementary immersions ($j$ open, $i$ closed)
of separated schemes of finite type over $\BC$. \\[0.1cm]
(a)~If $\bX$ is proper, then the result of applying $H^*(\bX(\BC),\argdot)$ 
to the boun\-dary triangle 
\[
j_! \BQ_X(0) \longto j_* \BQ_X(0) \longto i_* i^* j_* \BQ_X(0)
\longto j_! \BQ_X(0)[1] \; ,
\]
does only depend on $X$. \\[0.1cm]
(b)~The morphisms 
\[
i_* i^! \BQ_{\bX}(0) \longto i_* i^* \BQ_{\bX}(0)
\]
and
\[
i_* i^* \BQ_{\bX}(0) \longto i_* i^* j_* \BQ_X(0)
\]
induced by the respective adjunctions fit into a canonical exact triangle
\[
i_* i^! \BQ_{\bX}(0) \longto i_* i^* \BQ_{\bX}(0) \longto
i_* i^* j_* \BQ_X(0) \longto i_* i^! \BQ_{\bX}(0)[1] \; .
\]
It is the third column of a diagram of exact triangles
\[
\vcenter{\xymatrix@R-10pt{
        0 \ar[d] \ar[r] &
        i_* i^! \BQ_{\bX}(0) \ar[d] \ar@{=}[r] &
        i_* i^! \BQ_{\bX}(0) \ar[d] \ar[r]  &
        0 \ar[d] \\
        j_! \BQ_X(0) \ar@{=}[d] \ar[r] &
        \BQ_{\bX}(0) \ar[d] \ar[r] &
        i_* i^* \BQ_{\bX}(0) \ar[d] \ar[r] &
        j_! \BQ_X(0)[1] \ar@{=}[d] \\
        j_! \BQ_X(0) \ar[d] \ar[r] &  
        j_* \BQ_X(0) \ar[d] \ar[r] &
        i_* i^* j_* \BQ_X(0) \ar[d] \ar[r] &
        j_! \BQ_X(0)[1] \ar[d] \\
        0 \ar[r] &
        i_* i^! \BQ_{\bX}(0)[1] \ar@{=}[r] &
        i_* i^! \BQ_{\bX}(0)[1] \ar[r]  &
        0 
\\}}
\]
whose second and third rows are the localization and boundary triangles,
and whose second column is the co-localization triangle. \\[0.1cm]
(c)~If $\bX$ is smooth of constant dimension $d$, then there is a 
canonical isomorphism
\[
i_* i^! \BQ_{\bX}(0) \cong \BD_{\bX} \bigl( i_* i^* \BQ_{\bX}(d)[2d] \bigr) 
\] 
($\BD_{\bX} :=$ duality for Hodge modules on $\bX$ \cite[(4.1.5)]{Sa}). 
\end{Thm}

\begin{Proof}
(a) is a consequence of proper base change \cite[(4.4.3)]{Sa}. As for (b),
let us define the morphism
\[
i_* i^* j_* \BQ_X(0) \longto i_* i^! \BQ_{\bX}(0)[1]
\]
as $i_* i^*$ of the morphism 
\[
j_* \BQ_X(0) \longto i_* i^! \BQ_{\bX}(0)[1]
\]
occurring in the co-localization triangle. 
Together with the morphisms defined before, it yields the diagram
of the statement. Exactness of its third column is then 
a consequence of exactness of the first and second column.  
Finally, part~(c) results from duality \cite[(4.3.5)]{Sa}.
\end{Proof}

Fix $X \in Sch/\BC$.
In the sequel, the boundary triangle 
\[
j_! \BQ_X(0) \longto j_* \BQ_X(0) \longto i_* i^* j_* \BQ_X(0)
\longto j_! \BQ_X(0)[1] \; ,
\]
will always be assumed
to be formed using a compactification $j: X \into \bX$ of $X$,
with complement $i: D \into \bX$.
Theorem~\ref{1F} gives the precise relation between the boundary triangle
on the one hand, and the localization and co-localization triangles
on the other hand. While boundary cohomology, i.e., cohomology
of $i_* i^* j_* \BQ_X(0)$ does 
not depend on $\bX$, cohomology of the two other terms 
of the triangle 
\[
i_* i^! \BQ_{\bX}(0) \longto i_* i^* \BQ_{\bX}(0) \longto
i_* i^* j_* \BQ_X(0) \longto i_* i^! \BQ_{\bX}(0)[1] \; .
\]
from Theorem~\ref{1F}~(b) in general does 
(unless $X$ is itself proper). Saito's formalism allows to put restrictions
on the Hodge structures potentially occurring
as such cohomology groups.

\begin{Thm} \label{1G}
Let $n$ be an integer. Assume $\bX$ to be proper and smooth
(hence $X$ is smooth). \\[0.1cm]
(a)~The Hodge structure $H^n(\bX,i_* i^* \BQ_{\bX}(0))$ is of weights 
at most $n$. \\[0.1cm]
(b)~The Hodge structure $H^n(\bX,i_* i^! \BQ_{\bX}(0)[1])$ is of weights 
at least $n+1$.
\end{Thm}

\begin{Proof}
The scheme $\bX$ being proper, $H^n$ maps complexes of Hodge modules
of weights $\le 0$ to Hodge structures of weights $\le n$, and
complexes of Hodge modules
of weights $\ge 1$ to Hodge structures of weights $\ge n+1$. We thus need
to show that $i_* i^* \BQ_{\bX}(0)$ is of weights  $\le 0$, and 
$i_* i^! \BQ_{\bX}(0)$ of weights $\ge 0$. 

The scheme $\bX$ being smooth, $\BQ_{\bX}(0)$ is pure of weight $0$.
Therefore \cite[(4.5.2)]{Sa}, $i^* \BQ_{\bX}(0)$ is of weights  $\le 0$, and 
$i^! \BQ_{\bX}(0)$ of weights $\ge 0$, and the same remains true after
application of the functor $i_*$.
\end{Proof}

In ``triangulated'' language, Theorem~\ref{1G} says that the objects
\[
\bar{\pi}_* \bigl( i_* i^* \BQ_{\bX}(0) \bigr) \quad \text{and} \quad
\bar{\pi}_* \bigl( i_* i^! \BQ_{\bX}(0)[1] \bigr) 
\]
($\bar{\pi} :=$ the structure morphism of $\bX$)
of the derived category of Hodge structures are of weights $\le 0$
and $\ge 1$, respectively, when the compactification $\bX$ is smooth. 

\begin{Cor} \label{1H}
Let $X \in Sm/\BC$, and
\[
\ldots \longto A^n \longto \partial H^n (X (\BC),\BQ) \longto B^n
\]
\[
\longto A^{n+1} \longto \ldots \quad\quad\quad\quad\quad\quad
\]
a long exact sequence of mixed $\BQ$-Hodge structures. This sequence
is the result of applying $H^*(\bX(\BC),\argdot)$ 
to the triangle 
\[
i_* i^! \BQ_{\bX}(0) \longto i_* i^* \BQ_{\bX}(0) \longto
i_* i^* j_* \BQ_X(0) \longto i_* i^! \BQ_{\bX}(0)[1] \; ,
\]
for a suitable smooth compactification $j : X \into \bX$, only if
$A^n$ is of weights at most $n$, and $B^n$ is of weigths at least $n+1$,
for all $n \in \BZ$. 
\end{Cor}

Theorems~\ref{1F} and \ref{1G}, and Corollary~\ref{1H}
admit motivic analogues
(Theorems~\ref{1N} and \ref{1O}, and Corollary~\ref{1P} below), 
which we shall develop now.
To do so, it is necessary to use
the right notion of weights on triangulated categories.
Let us recall the following definitions and results of Bondarko
\cite{Bo}.

\begin{Def} \label{1I}
Let $\CC$ be a triangulated category. A \emph{weight structure on $\CC$}
is a pair $w = (\CC_{w \le 0} , \CC_{w \ge 0})$ of full 
sub-categories of $\CC$, such that, putting
\[
\CC_{w \le n} := \CC_{w \le 0}[n] \quad , \quad
\CC_{w \ge n} := \CC_{w \ge 0}[n] \quad \forall \; n \in \BZ \; ,
\]
the following conditions are satisfied.
\begin{enumerate}
\item[(1)] The categories
$\CC_{w \le 0}$ and $\CC_{w \ge 0}$ are 
Karoubi-closed: for any object $M$ of $\CC_{w \le 0}$ or
$\CC_{w \ge 0}$, any direct summand of $M$ formed in $\CC$
is an object of $\CC_{w \le 0}$ or
$\CC_{w \ge 0}$, respectively.
\item[(2)] (Semi-invariance with respect to shifts.)
We have the inclusions
\[
\CC_{w \le 0} \subset \CC_{w \le 1} \quad , \quad
\CC_{w \ge 0} \supset \CC_{w \ge 1}
\]
of full sub-categories of $\CC$.
\item[(3)] (Orthogonality.)
For any pair of objects $M \in \CC_{w \le 0}$ and $N \in \CC_{w \ge 1}$,
we have
\[
\Hom_{\CC}(M,N) = 0 \; .
\]
\item[(4)] (Weight filtration.)
For any object $M \in \CC$, there exists an exact triangle
\[
A \longto M \longto B \longto A[1]
\]
in $\CC$, such that $A \in \CC_{w \le 0}$ and $B \in \CC_{w \ge 1}$.
\end{enumerate}
\end{Def}

By condition~\ref{1I}~(2),
\[
\CC_{w \le n} \subset \CC_{w \le 0}
\]
for negative $n$, and
\[
\CC_{w \ge n} \subset \CC_{w \ge 0}
\]
for positive $n$. There are obvious analogues of the other conditions
for all the categories $\CC_{w \le n}$ and $\CC_{w \ge n}$. In particular,
they are all Karoubi-closed, and any object $M \in \CC$
is part of an exact triangle
\[
A \longto M \longto B \longto A[1]
\]
in $\CC$, such that $A \in \CC_{w \le n}$ and $B \in \CC_{w \ge n+1}$.
By a slight generalization of the terminology introduced in  
condition~\ref{1I}~(4), we shall refer to any such exact triangle
as a weight filtration of $M$.

\begin{Rem} \label{1K}
(a)~Our convention concerning the sign of the weight is actually opposite 
to the one from \cite[Def.~1.1.1]{Bo}, i.e., we exchanged the
roles of $\CC_{w \le 0}$ and $\CC_{w \ge 0}$. \\[0.1cm]
(b)~Note that in condition~\ref{1I}~(4), ``the'' weight filtration is not
assumed to be unique. 
\end{Rem}

As observed by Bondarko,
weight structures are relevant to motives thanks to the following
result.

\begin{Thm} \label{1L}
Let $F$ be a commutative flat $\BZ$-algebra,
and assume $k$ to admit resolution of singularities. \\[0.1cm]
(a)~There is a canonical weight structure on the category $\DeffgM_F$.
It is uniquely characterized by the requirement that its heart equal 
$\CHeffM_F$. \\[0.1cm]
(b)~There is a canonical weight structure on the category $\DgM_F$,
extending the weight structure from (a).
It is uniquely characterized by the requirement that its heart equal 
$\CHM_F$. \\[0.1cm]
(c)~Statements (a) and (b) hold without assuming resolution of singularities  
provided $F$ is a $\BQ$-algebra.
\end{Thm}

\begin{Proof}
For $F = \BZ$ and $k$ of characteristic zero,
this is the content of \cite[Sect.~6.5 and 6.6]{Bo}.
For the modifications of the proof in the remaining cases,
see \cite[Thm.~1.13]{W3}.
\end{Proof}

The following result is formally implied by Theorem~\ref{1L},
and the fundamental properties of the category $\DgM_F$,
notably \emph{localization} and 
\emph{duality} \cite[Prop.~4.1.5 and Thm.~4.3.7]{V}.
For details of the proof, we refer to \cite[Cor.~1.14]{W3} 
(\cite[Thm.~6.2.1~(1) and (2)]{Bo1} if $k$ is of characteristic zero).

\begin{Cor} \label{1M}
Let $X \in Sch/k$ be of (Krull) dimension $d$. 
Assume $k$ to admit resolution of singularities. \\[0.1cm]
(a)~The motive with compact support $\Mcgm(X)$ lies in 
\[
\DeffgM_{w \ge 0} \cap \DeffgM_{w \le d} \; . 
\]
(b)~If $X \in Sm/k$,
then the motive $M(X)$ 
lies in 
\[
\DeffgM_{w \ge -d} \cap \DeffgM_{w \le 0} \; .
\]
\end{Cor}

Fix $X \in Sch/k$. The motivic analogue of (the complex computing)
boundary cohomology (for $k= \BC$) is given by the 
\emph{boundary motive} $\dMgm(X)$ of $X$ \cite[Def.~2.1]{W1}.
The analogue of Theorem~\ref{1F} reads as follows;
there as in the sequel, we shall denote by $M^*$ the dual of a 
geometrical motive $M$ \cite[Thm.~4.3.7]{V}.

\begin{Thm} \label{1N}
(a)~There is a canonical exact boundary triangle
\[
\dMgm(X) \longto M(X) \longto \Mcgm(X) \longto \dMgm(X)[1]
\]
in $\DeffgM$. \\[0.1cm]
(b)~Assume $k$ to admit resolution of singularities. Let 
\[
\vcenter{\xymatrix@R-10pt{
        X \ar@{^{ (}->}[r]^-{j} &
        \bX \ar@{<-^{ )}}[r]^-{i} &
        D
\\}}
\]
be complementary immersions ($j$ open, $i$ closed)
of schemes in $Sch/k$. Assume $\bX$ to be proper. There is a
canonical morphism $\alpha: \dMgm(X) \to M(D)$. Define $M(\bX / X)$
as the relative motive of $\bX$ modulo $X$ \cite[Conv.~1.2]{W1},
and let $\beta: M(D) \to M(\bX / X)$ 
be induced by the morphism $i_*: M(D) \to M(\bX)$.
Then the morphisms $\alpha$ and $\beta$ fit into a canonical 
exact triangle
\[
M(\bX / X)[-1] \longto \dMgm(X) \stackrel{\alpha}{\longto} 
M(D) \stackrel{\beta}{\longto} M(\bX / X) \; .
\] 
It is the third column of a canonical
diagram of exact triangles
\[
\vcenter{\xymatrix@R-10pt{
        0 &
        M(\bX / X) \ar[l] &
        M(\bX / X) \ar@{=}[l]  &
        0 \ar[l] \\
        \Mcgm(X) \ar[u] &
        M(\bX) \ar[u] \ar[l] &
        M(D) \ar[u]_{\beta} \ar[l] &
        \Mcgm(X)[-1] \ar[u] \ar[l] \\
        \Mcgm(X) \ar@{=}[u]&  
        M(X) \ar[u] \ar[l] &
        \dMgm(X) \ar[u]_{\alpha} \ar[l] &
        \Mcgm(X)[-1] \ar@{=}[u] \ar[l] \\
        0 \ar[u] &
        M(\bX / X)[-1] \ar[u] \ar[l] &
        M(\bX / X)[-1] \ar[u] \ar@{=}[l] &
        0 \ar[l] \ar[u]   
\\}}
\]
whose second and third row are the localization \cite[Prop.~4.1.5]{V}
and boundary triangles,
and whose second column is 
cano\-nically associated to the relative motive $M(\bX / X)$. \\[0.1cm]
(c)~In the situation of (b), assume
$\bX$ to be proper and smooth of constant dimension $d$
(hence $X$ is smooth). 
There is a canonical morphism $\alpha^*: M(D)^*(d)[2d-1] \to \dMgm(X)$.
Let $\gamma:= i^*i_*: M(D) \to M(D)^*(d)[2d]$ be the composition
of $i_*: M(D) \to M(\bX)$ and its dual.
Then the morphisms $\alpha^*$, $\alpha$ (from (b)) 
and $\gamma$ form a canonical 
exact triangle
\[
M(D)^*(d)[2d-1] \stackrel{\alpha^*}{\longto} \dMgm(X) \stackrel{\alpha}{\longto} 
M(D) \stackrel{\gamma}{\longto} M(D)^*(d)[2d] \; .
\] 
It is the third column of a second canonical
diagram of exact triangles
\[
\vcenter{\xymatrix@R-10pt{
        0 &
        M(D)^*(d)[2d] \ar[l] &
        M(D)^*(d)[2d] \ar@{=}[l]  &
        0 \ar[l] \\
        \Mcgm(X) \ar[u] &
        M(\bX) \ar[u] \ar[l] &
        M(D) \ar[u]_{\gamma} \ar[l] &
        \Mcgm(X)[-1] \ar[u] \ar[l] \\
        \Mcgm(X) \ar@{=}[u]&  
        M(X) \ar[u] \ar[l] &
        \dMgm(X) \ar[u]_{\alpha} \ar[l] &
        \Mcgm(X)[-1] \ar@{=}[u] \ar[l] \\
        0 \ar[u] &
        M(D)^*(d)[2d-1] \ar[u] \ar[l] &
        M(D)^*(d)[2d-1] \ar[u]_{\alpha*} \ar@{=}[l] &
        0 \ar[l] \ar[u]   
\\}}
\]
whose second and third row are the localization 
and boundary triangles,
and whose second column is dual, up to a twist by $(d)$ 
and a shift by $[2d]$, to the localization triangle
(it will be referred to as the co-localization triangle). 
This diagram is isomorphic to the diagram from (b).
\end{Thm}

\begin{Proof}
For (a), let us briefly recall the definition of $\dMgm(X)$.
First, according to \cite[pp.~223, 224]{V},
a monomorphism of \emph{Nisnevich sheaves with transfers}
$\iota_X: L(X) \into L^c(X)$ is associated to $X$: the sheaf $L(X)$
is formed using \emph{finite correspondences}, and $L^c(X)$
is formed using \emph{quasi-finite correspondences}.
Next \cite[pp.~207, 208]{V}, there is a functor $\RC$ associating to a 
Nisnevich sheaf with transfers its \emph{singular simplicial complex}. 
Voevodsky goes on to define the
motive $M (X)$ as $\RC (L(X))$,
and the motive with compact support $\Mcgm (X)$ as
$\RC (L^c(X))$. Set 
\[
\dMgm (X) := \RC (\coker \iota_X) [-1] 
\]
\cite[Def.~2.1]{W1}. 
Claim~(a) is then a direct consequence of this definition. 
As for (b), we refer to \cite[Prop.~2.4]{W1}.
It remains to show part~(c).
The morphism $\alpha^*$ is defined as the dual of $\alpha$,
twisted by $(d)$ and shifted by $[2d-1]$
\[
M(D)^*(d)[2d-1] \longto \dMgm(X)^*(d)[2d-1] \; ,
\]
followed by the auto-duality 
isomorphism 
\[
\dMgm(X)^*(d)[2d-1] \isoto \dMgm(X)
\]
\cite[Thm.~6.1]{W1}.
Note that by duality \cite[Thm.~4.3.7~3]{V},
\[
M(\bX)^*(d)[2d] \isoto M(\bX)   \quad \text \quad  
M(X)^*(d)[2d] \isoto \Mcgm(X) 
\]
canonically, 
and under these identifications, the dual of the
canonical morphism $M(\bX) \to \Mcgm(X)$
occurring in the localization triangle
equals the canonical morphism $M(X) \to M(\bX)$
occurring in the co-localization triangle.
It remains to show that the composition
\[
M(D)^*(d)[2d-1] \stackrel{\alpha^*}{\longto} \dMgm(X) \longto M(X)
\]
equals the morphism
\[
M(D)^*(d)[2d-1] \longto M(X)
\]
in the co-localization sequence. But this identity can be checked after
applying duality. Note that the
boundary triangle is auto-dual \cite[Thm.~6.1]{W1}.
Therefore, the dual of the above composition equals
\[
\Mcgm(X) \longto \dMgm(X) \stackrel{\alpha}{\longto} M(D) \; ,
\]
which in turn equals the morphism
\[
\Mcgm(X) \longto M(D)
\]
in the localization sequence. 
\end{Proof}

Recall that motives \emph{\`a la} Voevodsky
behave homologically; this is why the sense of the 
arrows is inversed when compared to 
cohomology.

\begin{Rem}
If $\bX$ is proper and smooth of constant dimension,
there should be a \emph{canonical} choice of isomorphism
between the two canonical diagrams from Theorem~\ref{1N}~(b) and (c). 
If $D$ is (proper and) smooth, then such a choice is induced
by \emph{purity} \cite[Prop.~3.5.4]{V}.
\end{Rem}

Here is the motivic analogue of Theorem~\ref{1G};
it follows directly from Corollary~\ref{1M}.

\begin{Thm} \label{1O}
Assume $k$ to admit resolution of singularities. Let 
\[
\vcenter{\xymatrix@R-10pt{
        X \ar@{^{ (}->}[r]^-{j} &
        \bX \ar@{<-^{ )}}[r]^-{i} &
        D
\\}}
\]
be complementary immersions ($j$ open, $i$ closed)
of schemes in $Sch/k$. Assume $\bX$ to be proper and smooth
(hence $X$ is smooth). \\[0.1cm]
(a)~The motive $M(D)$ lies in $\DeffgM_{w \ge 0}$. \\[0.1cm]
(b)~The motive $M(D)^*(d)[2d-1]$ lies in $\DeffgM_{w \le -1}$. 
\end{Thm}

In particular, the exact triangle 
\[
M(D)^*(d)[2d-1] \stackrel{\alpha^*}{\longto} \dMgm(X) \stackrel{\alpha}{\longto} 
M(D) \stackrel{\gamma}{\longto} M(D)^*(d)[2d] \; .
\] 
from Theorem~\ref{1N}~(c) is then a weight filtration of $\dMgm(X)$.

\begin{Cor} \label{1P}
Assume $k$ to admit resolution of singularities.
Let 
\[
A \longto \dMgm(X) \longto B \longto A[1]
\]
be an exact triangle in $\DeffgM$, for $X \in Sm/k$. This triangle
is isomorphic to the triangle 
\[
M(D)^*(d)[2d-1] \stackrel{\alpha^*}{\longto} \dMgm(X) \stackrel{\alpha}{\longto} 
M(D) \stackrel{\gamma}{\longto} M(D)^*(d)[2d] \; .
\] 
for a suitable smooth compactification $j : X \into \bX$, only if
it is a weight filtration of $\dMgm(X)$:
$A \in \DeffgM_{w \le -1}$ and $B \in \DeffgM_{w \ge 0}$. 
\end{Cor}

Altogether, for fixed $X \in Sm/k$, we get a functor
from the category of smooth compactifications of $X$
to the category of weight filtrations of $\dMgm(X)$.
It turns out to be very instructive to see what
one gets when trying to invert this functor.

\begin{Constr} \label{1Q}
Assume $k$ to admit resolution of singularities.
Fix a weight filtration 
\[
\dMgm(X)_{\le -1} \stackrel{c_-}{\longto} \dMgm(X) 
\stackrel{c_+}{\longto} \dMgm(X)_{\ge 0} 
\stackrel{\delta}{\longto} \dMgm(X)_{\le -1}[1]
\]
of $\dMgm(X)$, for $X \in Sm/k$:
\[
\dMgm(X)_{\le -1} \in \DeffgM_{w \le -1} \;\; \text{and} \quad 
\dMgm(X)_{\ge 0} \in \DeffgM_{w \ge 0} \; .
\] 
Consider the boundary triangle
\[
\dMgm(X) \stackrel{v_-}{\longto} M(X) 
\stackrel{u}{\longto} \Mcgm(X) \stackrel{v_+}{\longto} \dMgm(X)[1] \; .
\] 
According to axiom TR4' of triangulated categories (see \cite[Sect.~1.1.6]{BBD}
for an equivalent formulation), the diagram of exact triangles
\[
\vcenter{\xymatrix@R-10pt{
        0 &
        \dMgm(X)_{\le -1}[1] \ar[l] &
        \dMgm(X)_{\le -1}[1] \ar@{=}[l]  &
        0 \ar[l] \\
        \Mcgm(X) \ar[u] &
         &
        \dMgm(X)_{\ge 0} \ar[u]_{\delta} &
        \Mcgm(X)[-1] \ar[u] \ar[l]_{c_+(v_+[-1])} \\
        \Mcgm(X) \ar@{=}[u]&  
        M(X) \ar[l]_-u &
        \dMgm(X) \ar[u]_{c_+} \ar[l]_-{v_-} &
        \Mcgm(X)[-1] \ar@{=}[u] \ar[l]_-{v_+[-1]} \\
        0 \ar[u] &
        \dMgm(X)_{\le -1} \ar[u]^{v_-c_-} \ar[l] &
        \dMgm(X)_{\le -1} \ar[u]_{c_-} \ar@{=}[l] &
        0 \ar[l] \ar[u]   
\\}}
\]
can be completed to give
\[
\vcenter{\xymatrix@R-10pt{
        0 &
        \dMgm(X)_{\le -1}[1] \ar[l] &
        \dMgm(X)_{\le -1}[1] \ar@{=}[l]  &
        0 \ar[l] \\
        \Mcgm(X) \ar[u] &
        M_0 \ar[u]^{\delta_-} \ar[l]_-{i_0}  &
        \dMgm(X)_{\ge 0} \ar[u]_{\delta} \ar[l]_-{\delta_+} &
        \Mcgm(X)[-1] \ar[u] \ar[l]_{c_+(v_+[-1])} \\
        \Mcgm(X) \ar@{=}[u]&  
        M(X) \ar[u]^{\pi_0} \ar[l]_-u &
        \dMgm(X) \ar[u]_{c_+} \ar[l]_-{v_-} &
        \Mcgm(X)[-1] \ar@{=}[u] \ar[l]_-{v_+[-1]} \\
        0 \ar[u] &
        \dMgm(X)_{\le -1} \ar[u]^{v_-c_-} \ar[l] &
        \dMgm(X)_{\le -1} \ar[u]_{c_-} \ar@{=}[l] &
        0 \ar[l] \ar[u]   
\\}}
\]
Note that this completion necessitates \emph{choices} of $M_0$
and of factorizations $u = i_0 \pi_0$ and $\delta = \delta_- \delta_+$.
In general, the object $M_0$ is unique up to possibly non-unique
isomorphism; it is this problem that will be addressed in the last part
of this section. \\

For the moment, note that whatever choice we make, $M_0$ will 
be in the heart of our weight structure: indeed, the second row
of the diagram, together with Corollary~\ref{1M}~(a) shows
that $M_0 \in \DeffgM_{w \ge 0}$, and the second column,
together with Corollary~\ref{1M}~(b) shows that $M_0 \in \DeffgM_{w \le 0}$.
According to Theorem~\ref{1L}~(a), it is therefore an effective
Chow motive. Note that it comes equipped with a factorization
\[
M(X) \stackrel{\pi_0}{\longto} M_0
\stackrel{i_0}{\longto} \Mcgm(X)
\]
of the canonical morphism $u: M(X) \to \Mcgm(X)$, and that the triangles
\[
\dMgm(X)_{\ge 0} \stackrel{\delta_+}{\longto} M_0
\stackrel{i_0}{\longto} \Mcgm(X) 
\stackrel{(c_+[1])v_+}{\longto} \dMgm(X)_{\ge 0}[1]
\]
and
\[
\dMgm(X)_{\le -1} \stackrel{v_-c_-}{\longto} M(X)
\stackrel{\pi_0}{\longto} M_0
\stackrel{\delta_-}{\longto} \dMgm(X)_{\le -1}[1]
\]
are weight filtrations of $\Mcgm(X)$ and of $M(X)$, respectively.
\end{Constr}

Let us summarize the discussion.

\begin{Thm} \label{1R}
Assume $k$ to admit resolution of singularities, and fix $X \in Sm/k$.
The map
\[
\big\{ \bigl( \dMgm(X)_{\le -1} , \dMgm(X)_{\ge 0} \bigr) \big\} / \cong \; 
\longto 
\big\{ \bigl( M_0 , \pi_0 , i_0 \bigr) \big\} / \cong
\]
from the preceding construction
is a bijection between \\[0.1cm]
(1)~the isomorphism classes of weight filtrations of 
the boundary motive $\dMgm(X)$, \\[0.1cm]
(2)~the isomorphism classes of effective Chow motives
$M_0$, together with a factorization 
\[
M(X) \stackrel{\pi_0}{\longto} M_0
\stackrel{i_0}{\longto} \Mcgm(X)
\]
of the canonical morphism $u: M(X) \to \Mcgm(X)$, such that 
both $i_0$ and $\pi_0$ can be completed to give weight filtrations of 
$\Mcgm(X)$ and of $M(X)$, respectively.
\end{Thm}

There are obvious $F$-linear versions of Theorem~\ref{1R},
for any commutative flat $\BZ$-algebra $F$.
Recall that we started off with special choices of Chow motives
factorizing $u$, namely the motives of smooth compactifications
of $X$. But Theorem~\ref{1R} should yield more
general Chow motives $M_0$. For example, one might hope for
the motivic version of intersection cohomology of a singular
compactification of $X$ to occur. For surfaces, this will 
be spelled out in Section~\ref{2}. \\

Note that we are forced to pass to the level of
isomorphism classes because of the choices made in Construction~\ref{1Q}.
One important problem caused by this is the lack of functoriality.
Thus, an endomorphism of a given weight filtration
\[
\dMgm(X)_{\le -1} \longto \dMgm(X) 
\longto \dMgm(X)_{\ge 0} 
\longto \dMgm(X)_{\le -1}[1]
\]
will in general not yield an endomorphism 
of any of the Chow motives $M_0$ representating 
the associated isomorphism class. 

\begin{Prin} \label{1Ra}
In order to obtain functoriality, Construction~\ref{1Q}
needs to be \emph{rigidified}. 
\end{Prin}

It turns out that an \emph{ad hoc} 
geometrical method
suffices to achieve rigidification in the setting of surfaces
(see Section~\ref{2}). Let us finish this section by describing
another method (namely, that of \cite{W3}), 
based again on the formalism of weights.
It will be illustrated in the setting of 
self-products of the universal elliptic curve over 
a modular curve (see Section~\ref{3}). \\

\begin{Rem} \label{1S}
Getting back to Construction~\ref{1Q}, and starting again
with a weight filtration
\[
\dMgm(X)_{\le -1} \stackrel{c_-}{\longto} \dMgm(X) 
\stackrel{c_+}{\longto} \dMgm(X)_{\ge 0} 
\stackrel{\delta}{\longto} \dMgm(X)_{\le -1}[1] \; ,
\]
let us see what obstacles
there are for the triple $(M_0,\pi_0,i_0)$ to be unique up to \emph{unique}
isomorphism. Note that $(M_0,\pi_0)$ is a cone of
\[
v_-c_-: \dMgm(X)_{\le -1} \longto \dMgm(X) \longto M(X) \; .
\]
Any other choice of cone would map isomorphically to $(M_0,\pi_0)$,
the isomorphism in question being unique up to the image of
an element in
\[
\Hom_{\DeffgM} \bigl( \dMgm(X)_{\le -1}[1] , M_0 \bigr) \; .
\]
In general, the object $\dMgm(X)_{\le -1}$ belonging to $\DeffgM_{w \le -1}$,
hence 
\[
\dMgm(X)_{\le -1}[1] \in \DeffgM_{w \le 0} \; ,
\]
there is no way of preventing such elements from being non-zero.
However, if   
\[
\dMgm(X)_{\le -1} \in \DeffgM_{w \le -2} \subset \DeffgM_{w \le -1} \; ,
\]
then 
\[
\Hom_{\DeffgM} \bigl( \dMgm(X)_{\le -1}[1] , M_0 \bigr) = 0
\]
by orthogonality~\ref{1I}~(3) (recall that $M_0$ is of weight zero).
Thus, under this hypothesis, the pair $(M_0,\pi_0)$ is rigid.
As for $i_0$, the same type of reasoning shows unicity provided that
\[
\dMgm(X)_{\ge 0} \in \DeffgM_{w \ge 1} \subset \DeffgM_{w \ge 0} \; .
\]
\end{Rem}

We are thus led naturally to make the following definition 
\cite[Def.~1.6 and 1.10]{W3}.

\begin{Def} \label{1T}
Let $M \in \DeffgM$, and $m \le n$ two integers (which may be identical). 
A \emph{weight filtration of $M$ avoiding weights $m,m+1,\ldots,n-1,n$}
is an exact triangle
\[ 
M_{\le {m-1}} \longto M \longto M_{\ge {n+1}} \longto M_{\le m-1}[1] \; ,
\]
with $M_{\le {m-1}} \in \DeffgM_{w \le m-1}$
and $M_{\ge {n+1}} \in \DeffgM_{w \ge {n+1}} \ $. 
If such a weight filtration exists, then we say that \emph{$M$ is without 
weights $m,\ldots,n$}.
\end{Def}

Weight filtrations avoiding weights $m,\ldots,n$ behave
functorially \cite[Prop.~1.7]{W3}. In particular, if
$M \in \DeffgM$ is without weights $m,\ldots,n$, then its weight filtration
avoiding weights $m,\ldots,n$ is unique up to unique ismorphism.
Remark~\ref{1S} therefore shows that we can rigidify 
Construction~\ref{1Q} provided that the boundary motive $\dMgm(X)$
is without weights $-1$ and $0$. 

\begin{Prob} \label{1U} 
The boundary motive $\dMgm(X)$ of $Sm/k$
is never without weights $-1$ and $0$
--- unless it is altogether trivial. \\

Here is a heuristic reason, using the weights occurring in boundary
cohomology over $k = \BC$: 
to say that $\dMgm(X)$ is not trivial implies that $X$ is
not proper. On the one hand, the cokernel of 
\[
H^0_c (X (\BC),\BQ) \longto H^0 (X (\BC),\BQ) \cong \BQ(0)^r
\]
is then non-trivial. On the other hand, it injects into 
$\partial H^0 (X (\BC),\BQ)$.
\end{Prob}

Therefore, the approach of ``avoiding weights'' cannot work 
on the whole of the boundary motive. We need to restrict to
direct factors. Fix a commutative flat $\BZ$-algebra $F$, and let  
\[
(\ast) \quad\quad \dMgm(X) \longto M(X) \longto \Mcgm(X) \longto \dMgm(X)[1]
\]
be the boundary triangle associated to a fixed object $X$ of $Sm/k$,
viewed as a triangle in $\DeffgM_F$.
Fix an idempotent endomorphism $e$ of the triangle $(\ast)$, that is,
fix idempotent endomorphisms of each of $M(X)$, $\Mcgm(X)$ and $\dMgm(X)$, 
viewed as objects of $\DeffgM_F$, which 
yield an endomorphism of the triangle.
Denote by $M(X)^e$, $\Mcgm(X)^e$ and $\dMgm(X)^e$ the images of $e$
on $M(X)$, $\Mcgm(X)$ and $\dMgm(X)$, respectively, 
and consider the canonical morphism
$u: M(X)^e \to \Mcgm(X)^e$. By Corollary~\ref{1M}
and condition \ref{1I}~(1), the object $M(X)^e$
belongs to $\DeffgM_{F,w \le 0}$, and 
$\Mcgm(X)^e$ to $\DeffgM_{F,w \ge 0} \ $. 
In this situation, Construction~\ref{1Q} yields the following.

\begin{Var} \label{1V}
The map
\[
\big\{ \bigl( \dMgm(X)^e_{\le -1} , \dMgm(X)^e_{\ge 0} \bigr) \big\} / \cong \;
\longto 
\big\{ \bigl( M_0 , \pi_0 , i_0 \bigr) \big\} / \cong
\]
is a bijection between \\[0.1cm]
(1)~the isomorphism classes of weight filtrations of $\dMgm(X)^e$, \\[0.1cm]
(2)~the isomorphism classes of effective Chow motives
$M_0 \in \CHeffM_F$, together with a factorization 
\[
M(X)^e \stackrel{\pi_0}{\longto} M_0
\stackrel{i_0}{\longto} \Mcgm(X)^e
\]
of the canonical morphism $u: M(X)^e \to \Mcgm(X)^e$, such that 
both $i_0$ and $\pi_0$ can be completed to give weight filtrations of 
$\Mcgm(X)^e$ and of $M(X)^e$, respectively.
\end{Var}

\begin{Comp} \label{1W}
Let $e$ denote an idempotent endomorphism of the boundary triangle
$(\ast)$ as above, and
assume that $\dMgm(X)^e$
is without weights $-1$ and $0$. 
Then the isomorphism class $(M_0,\pi_0,i_0)$ 
associated to the weight filtration of $\dMgm(X)^e$
avoiding weights $-1$ and $0$ essentially contains
one single object, which is unique up to unique isomorphism.
\end{Comp}

This is the principle exploited in \cite[Sect.~4]{W3}. 
Due to the behaviour of its realizations,
the object $M_0$ is referred to as the 
\emph{$e$-part of the interior motive} of $X$. It has very strong
functoriality properties. They will be illustrated in our Section~\ref{3},
where we shall establish equivariance under the Hecke algebra of $M_0$
in a special geometrical context. 
Before that, we need to address two very concrete questions: \\[0.1cm]
(I)~How does one get endomorphisms of the boundary triangle
\[
(\ast) \quad\quad 
\dMgm(X) \longto M(X) \longto \Mcgm(X) \longto \dMgm(X)[1] \quad\quad ?
\]
(II)~How can one show that a given such endomorphism is idempotent? \\[0.1cm]
The following two sections attempt to answer these questions,
at least partially.

%%% Local Variables:
%%% mode: latex
%%% TeX-master: "head"
%%% End:

\bigskip

%%%%%%%%%%%%%%%%%%%%%%%%%%%%%%%%%%%%%%%%%%%%%%%%%%%%%%%%%%%%%%%%%%%%%%%
%
%  Section 2b
%
%%%%%%%%%%%%%%%%%%%%%%%%%%%%%%%%%%%%%%%%%%%%%%%%%%%%%%%%%%%%%%%%%%%%%%%

\section{Relative motives and functoriality of the boundary motive}
\label{2b}

%%%%%%%%%%%%%%%%%%%%%%%%%%%%%%%%

%%%%%%%%%%%%%%%%%%%%%%%%%%%%%%%%

In this and the next section, the base field $k$ is assumed to admit
strict resolution of singularities.
For $X \in Sm/k$,
the algebra of \emph{finite correspondences} 
$c(X,X)$ acts on $M(X)$ \cite[p.~190]{V}. 
In order to apply the constructions from Section~\ref{1}, we need
to construct endomorphisms of the whole boundary triangle
\[
(\ast) \quad\quad
\dMgm(X) \longto M(X) \longto \Mcgm(X) \longto \dMgm(X)[1] \; .
\]
One of the aims of this section
is to show that the theory of \emph{relative
motives} provides a source of such endomorphisms. 
This result is a special feature
of an analysis of the functorial behavior
of the exact triangle $(\ast)$ under morphisms of relative motives
(Theorems~\ref{2bB} and \ref{2bBa}, Corollary~\ref{2bK}). 
The main application (Example~\ref{2bL})
concerns endomorphisms of $(*)$ ``of Hecke type''. \\

Let us fix a base scheme $S \in Sm/k$.  Recall that
by definition, objects of $Sm/k$ are separated
over $k$. Thus, for any two schemes $X$ and $Y$ over $S$, the natural morphism
\[
X \times_S Y \longto X \times_k Y
\]
is a closed immersion. Therefore, cycles on $X \times_S Y$
can and will be considered as cycles on $X \times_k Y$.
Denote by $Sm/S$ the category of separated smooth schemes of finite 
type over $S$, by $PropSm/S \subset Sm/S$ the full sub-category
of objects which are proper and smooth over $S$, and by
$ProjSm/S \subset Sm/S$ the full sub-category
of projective, smooth $S$-schemes. 

\begin{Def} \label{2bA}
Let $X,Y \in Sm/S$. 
Denote by $c_S(X,Y)$ the subgroup of $c(X,Y)$ of correspondences
whose support is contained in $X \times_S Y$. 
\end{Def}

The group $c_S(X,Y)$ is at the base of the theory 
of \emph{(effective) geometrical motives over $S$},
as defined and developed (for arbitrary regular Noetherian bases $S$) 
in \cite{Deg1,Deg2}. Note that any cycle $\FZ$ in $c_S(X,Y)$
gives rise to a morphism from $M(X)$ to $M(Y)$, which
we shall denote by $M(\FZ)$. 
Recall from \cite[Sect.~1.3, 1.6]{DM} the definition of
the categories of \emph{smooth (effective) Chow motives over $S$}; 
note that the approach of \loccit \ does not necessitate 
passage to $\BQ$-coefficients, and that one may choose to perform
the construction using schemes in $PropSm/S$ instead
of just schemes in $ProjSm/S$.
Denote by $\CHSeffM$ and $\CHSM$ the respective opposites
of these cate\-gories.
Note that for $X, Y \in PropSm/S$ and $\FZ \in c_S(X,Y)$,
the class of $\FZ$ in the Chow group $\ch^*(X \times_S Y)$
of cycles modulo rational equivalence
lies in the right degree, and therefore defines a morphism
from the relative Chow motive $h(X/S)$ of $X$ to the
relative Chow motive $h(Y/S)$.
Our aim is to prove the following.

\begin{Thm} \label{2bB}
(a)~There is a canonical additive covariant functor,
denoted $( \dMgm , M , \Mcgm ) = ( \dMgm , M , \Mcgm )_S$, from $\CHSM$
to the cate\-gory of exact triangles in $\DgM$.
On objects, it is characterized by the following properties:
\begin{enumerate}
\item[(a1)] for $X \in PropSm/S$, 
the functor $( \dMgm , M , \Mcgm )$ maps $h(X/S)$ to the 
boundary triangle
\[
(\ast)_X \quad\quad
\dMgm(X) \longto M(X) \longto \Mcgm(X) \longto \dMgm(X)[1] \; ,
\]
\item[(a2)] the functor $( \dMgm , M , \Mcgm )$ 
is compatible with Tate twists.
\end{enumerate}
On morphisms, the functor $( \dMgm , M , \Mcgm )$ 
maps the class of a cycle $\FZ \in c_S(X,Y)$
in $\ch^*(X \times_S Y)$, for $X,Y \in PropSm/S$, to 
a morphism $(\ast)_X \to (\ast)_Y$ whose $M$-component 
$M(X) \to M(Y)$ coincides
with $M(\FZ)$. \\[0.1cm]
(b)~There is a canonical additive contravariant functor
$( \dMgm , M , \Mcgm )^* = ( \dMgm , M , \Mcgm )^*_S$ from $\CHSM$
to the category of exact triangles in $\DgM$.
On objects, it is characterized by the following properties:
\begin{enumerate}
\item[(b1)] for an object $X \in PropSm/S$
which is of pure absolute dimension $d_X$, 
the functor $( \dMgm , M , \Mcgm )^*$ maps $h(X/S)$ to the triangle
\[
(\ast)^*_X := (\ast)_X(-d_X)[-2d_X] \; ,
\]
\item[(b2)] the functor $( \dMgm , M , \Mcgm )^*$ 
is anti-compatible with Tate twists.
\end{enumerate}
On morphisms, the functor $( \dMgm , M , \Mcgm )^*$ 
maps the class of a cycle $\FZ \in c_S(X,Y)$
in $\ch^*(X \times_S Y)$, for $X,Y \in PropSm/S$ of
pure absolute dimensions $d_X$ and $d_Y$, respectively, to 
a morphism $(\ast)^*_Y \to (\ast)^*_X$ 
whose $\Mcgm$-component coincides
with the dual of $M(\FZ)$. \\[0.1cm]
(c)~The functor $( \dMgm , M , \Mcgm )^*$ is canonically identified with the
composition of $( \dMgm , M , \Mcgm )$ and duality in $\DgM$.
\end{Thm}

Recall that by \cite[Thm.~4.3.7~3]{V},
the object $\Mcgm(X)$ is 
indeed dual to $M(X)(-d_X)[-2d_X]$. 
Note also \cite[Cor.~4.1.6]{V} that the functor 
from Theorem~\ref{2bB}~(a) maps the full sub-category $\CHSeffM$
to the full sub-category \cite[Thm.~4.3.1]{V} of exact triangles in $\DeffgM$.
Also recall that by convention, the Tate twist $(n)$ in $\CHSM$
corresponds to the (componentwise)
operation $M \mapsto M(n)[2n]$ in $\DgM$.
Thus, anti-compatibility of the functor $( \dMgm , M , \Mcgm )^*$
with Tate twists means that for any object $X$ of $\CHSM$, 
there is a functorial isomorphism
\[
\bigl( \dMgm , M , \Mcgm \bigr)^* (X(n)) \isoto
\bigl( \bigl( \dMgm , M , \Mcgm \bigr)^* (X) \bigr) (-n)[-2n] \; .
\]

\begin{Rem} \label{2bC}
As far as the $M$- and $\Mcgm$-components are concerned,
Theorem~\ref{2bB}, or at least its restriction to the full sub-category
$\CHSM_{proj}$ of $\CHSM$ 
generated by the motives of projective smooth $S$-schemes, 
is a consequence of the main results of \cite{Deg2},
especially \cite[Thm.~5.23]{Deg2}, together with the existence
of an adjoint pair $(L a_{S,\sharp}, a_S^*)$ of exact functors 
\cite[Ex.~4.12, Ex.~7.15]{CDeg} linking 
the category $\DSgM$ of geometrical motives over $S$ to $\DgM$
(here we let $a_S: S \to \Spec k$ denote the structure morphism of $S$). 
We should also mention that this approach would allow to avoid
the hypothesis on strict resolution of singularities.
However, the application of 
the results of \loccit \ to the functor $\dMgm$ is not obvious.
We are therefore forced to follow an alternative approach.
\end{Rem}

\begin{Rem}
The following
sheaf-theoretical phenomenon explains why one should expect
a statement like Theorem~\ref{2bB}. 
Writing $a = a_X$ for the structure
morphism $X \to \Spec k$, for $X \in Sch/k$, there is an exact
triangle of exact functors
\[
(+)_X \quad\quad a_! \longto a_* \longto a_* / a_! \longto a_![1]
\]
from the derived category $D^+(X)$ of complexes of \'etale sheaves on $X$
(say), bounded from below, to $D^+(\Spec k)$. 
Here, $a_*$ denotes the derived functor of the direct image, $a_!$
is its analogue ``with compact support'', and $a_* / a_!$
is a canonical choice of cone (which exists
since the category of compactifications of $X$
is filtered). The triangle $(+)_X$ is contravariantly functorial
with respect to proper morphisms. Up to a twist and a shift,
it is covariantly functorial with respect to proper smooth
morphisms. This shows that 
a suitable version of Theorem~\ref{2bB}~(a)
is likely to extend to the sub-category of   
$\DSgM$ generated by the relative motives of schemes 
which are (only) proper over $S$.
\end{Rem}

For any proper smooth morphism
$f: T \to S$ in the category $Sm/k$, denote by $f_\sharp: \CHTM \to \CHSM$
the canonical functor
induced by $h(X/T) \mapsto h(X/S)$, for any proper smooth
scheme $X$ over $T$ (hence, over $S$). For any morphism
$g: U \to S$ in $Sm/k$, denote by $g^*: \CHSM \to \CHUM$
the canonical tensor functor 
induced by $h(Y/S) \mapsto h(Y \times_S U/U)$, for any proper smooth
scheme $Y$ over $S$. 
When $g$ is proper and smooth, the functor $g_\sharp$ is left adjoint
to $g^*$. 
The following summarizes
the behaviour of 
$( \dMgm , M , \Mcgm )$ and $( \dMgm , M , \Mcgm )^*$
under change of the base $S$. 

\begin{Thm} \label{2bBa}
(a)~Let $f: T \to S$ be a proper smooth morphism in $Sm/k$.
There are canonical isomorphisms of additive functors
\[
\alpha_{f_{\sharp}} : ( \dMgm , M , \Mcgm )_S \circ f_\sharp \isoto 
( \dMgm , M , \Mcgm )_T 
\]
and
\[
\alpha_{f_{\sharp}}^* : ( \dMgm , M , \Mcgm )_T^* \isoto
( \dMgm , M , \Mcgm )_S^* \circ f_\sharp  
\]
on $\CHTM$.
The formation of both $\alpha_{f_{\sharp}}$ 
and $\alpha_{f_{\sharp}}^*$ is compatible with composition
of proper smooth morphisms in $Sm/k$. 
Under the identification of Theorem~\ref{2bB}~(c),
the equivalence $\alpha_{f_{\sharp}}^*$ corresponds to
the dual of the equivalence $\alpha_{f_{\sharp}}$. \\[0.1cm]
(b)~Let $g: U \to S$ be a proper smooth morphism in $Sm/k$.
Then there exists a canonical transformation of additive functors
\[
\beta_{g^*,\id_S} : ( \dMgm , M , \Mcgm )_U \circ g^* \longto
( \dMgm , M , \Mcgm )_S \; .
\]
The formation of $\beta_{g^*,\id_S}$ is compatible with composition
of proper smooth morphisms in $Sm/k$. \\[0.1cm]
(c)~The transformations $\alpha_{f_\sharp}$ and $\beta_{g^*,\id_S}$
commute in the following sense: let $f: T \to S$ and $g: U \to S$ be 
proper smooth morphisms in $Sm/k$. Consider the cartesian diagram
\[
\vcenter{\xymatrix@R-10pt{
        V  = T \times_S U\ar[d]_{g'} \ar[r]^-{f'} &
        U \ar[d]^g  \\
        T \ar[r]^f &
        S 
\\}}
\]
and the canonical identification of natural transformations
\[
f'_\sharp \circ g'^* = g^* \circ f_\sharp  
\]
of functors from $\CHTM$ to $\CHUM$.
Then the transformations
\[
\beta_{g'^*,\id_T} \circ (\alpha_{f'_\sharp} \circ g'^*) \quad , \quad
\alpha_{f_\sharp} \circ (\beta_{g^*,\id_S} \circ f_\sharp)
\]
of functors on $\CHTM$
\[
( \dMgm , M , \Mcgm )_U \circ g^* \circ f_\sharp \longto
( \dMgm , M , \Mcgm )_T
\]
coincide. \\[0.1cm]
(d)~Let $g: U \to S$ be a proper smooth morphism in $Sm/k$.
Then there exists a canonical transformation of additive functors
\[
\gamma_{\id_S,g^*} : ( \dMgm , M , \Mcgm )^*_S \longto 
( \dMgm , M , \Mcgm )^*_U \circ g^* \; .
\]
The formation of $\gamma_{\id_S,g^*}$ is compatible with composition
of proper smooth morphisms in $Sm/k$. 
Under the identification of Theorem~\ref{2bB}~(c),
the transformation $\gamma_{\id_S,g^*}$ corresponds to
the dual of the transformation $\beta_{g^*,\id_S}$.\\[0.1cm]
(e)~The transformations $\alpha_{f_\sharp}^*$ and $\gamma_{\id_S,g^*}$
commute in the following sense: let $f: T \to S$ and $g: U \to S$ be 
proper smooth morphisms in $Sm/k$.
Consider the cartesian diagram
\[
\vcenter{\xymatrix@R-10pt{
        V  = T \times_S U\ar[d]_{g'} \ar[r]^-{f'} &
        U \ar[d]^g  \\
        T \ar[r]^f &
        S 
\\}}
\]
Then the transformations
\[
(\gamma_{\id_S,g^*} \circ f_\sharp) \circ \alpha_{f_\sharp}^* \quad , \quad
(\alpha_{f'_\sharp}^* \circ g'^*) \circ \gamma_{\id_T,g'^*}  
\]
of functors on $\CHTM$
\[
( \dMgm , M , \Mcgm )_T^* \longto
( \dMgm , M , \Mcgm )_U^* \circ g^* \circ f_\sharp 
\]
coincide. 
\end{Thm}

\begin{Rem}
Sheaf-theoretical considerations show that parts (b)--(e)
of Theorem~\ref{2bBa}
should hold more generally for morphisms $g$ which are
(only) proper. While this could be shown to be indeed the case,
we chose to prove the statements only under the more restrictive
assumption on $g$: the proof then simplifies considerably since it
is possible to make use of the functor $g_\sharp \, $, which only exists
when $g$ is (proper and) smooth.
\end{Rem}

Let us prepare the proofs of Theorems~\ref{2bB} and \ref{2bBa}.
They are based on the following result.

\begin{Thm}[{\cite[Thm.~6.14, Rem.~6.15]{W1}}] \label{2bD}
Let $W \in Sm/k$ be of pure dimension $m$, and
$Z \subset W$ a closed sub-scheme
such that arbitrary intersections of
the irreducible components of $Z$ are smooth. Fix $n \in \BZ$. \\[0.1cm]
(a)~There is a canonical morphism
\[
cyc: \uh^0 \bigl( \zeq (W, m-n)_Z \bigr) (\Spec k) \longto  
\Hom_{\DeffgM} (M (W / Z), \BZ(n)[2n]) \; .
\]
(b)~The morphism $cyc$ is compatible with passage from
the pair $Z \subset W$ to $Z' \subset U$, for open
sub-schemes $U$ of $W$, and closed sub-schemes $Z'$ of $Z \cap U$
such that arbitrary intersections of
the irreducible components of $Z'$ are smooth. \\[0.1cm]
(c)~When $Z$ is empty, then  
\[
cyc: \uh^0 \bigl( \zeq (W, m-n) \bigr) (\Spec k) \longto  
\Hom_{\DeffgM} (M (W), \BZ(n)[2n]) 
\]
coincides with the morphism
from \cite[Cor.~4.2.5]{V}. In particular, it is then an isomorphism.
\end{Thm}

Some explanations are necessary. First,
by definition \cite[Def.~6.13]{W1}, the Nisnevich sheaf with transfers 
$\zeq (W , m-n)_Z$ associates to
$T \in Sm/k$ the group of those cycles in $\zeq (W , m-n)(T)$
\cite[p.~228]{V} having empty intersection with $T \times_k Z$.
In particular, the group $\zeq (W , m-n)_Z(\Spec k)$
equals the group of cycles on $W$ of dimension $m-n$,
whose support is disjoint from $Z$.
Recall then \cite[p.~207]{V} that the group
\[ 
\uh^0 \bigl( \zeq (W , m-n)_Z \bigr) (\Spec k) 
\]
is the quotient of $\zeq (W , m-n)_Z (\Spec k)$
by the image under the differential ``pull-back \emph{via} $1$ minus
pull-back \emph{via} $0$'' of $\zeq (W , m-n)_Z (\BA^1_k)$.
Finally the object $M (W / Z)$ denotes the relative motive
associated to the immersion of $Z$ into $W$ \cite[Def.~6.4]{W1}. 

\begin{Rem}
One may speculate about the validity of Theorem~\ref{2bD}
for arbitrary closed sub-schemes $Z$ of $W \in Sm/S$. 
While the author is optimistic about this possibility,
he notes that the tools developed in \cite{W1} to prove 
Theorem~\ref{2bD} require $Z$ to satisfy our more restrictive
hypotheses. It is for that reason that we are forced to
suppose $k$ to admit strict resolution of singularities.
\end{Rem}

Now note the following.

\begin{Prop} \label{2bE}
In the above situation, let in addition $V \subset W$ be a closed
sub-scheme in $Sm/S$, which is disjoint from $Z$.
Then the natural map
\[
\zeq (V , m-n) \longto \zeq (W , m-n)_Z 
\]
induces a morphism
\[
\ch_{m-n}(V) \longto \uh^0 \bigl( \zeq (W , m-n)_Z \bigr) (\Spec k) \; .
\]
\end{Prop}

\begin{Cor} \label{2bF}
Let $W \in Sm/k$ be of pure dimension $m$, 
$V, Z \subset W$ closed sub-schemes, and $n \in \BZ$.
Suppose that arbitrary intersections of
the irreducible components of $Z$ are smooth,
and that $V \cap Z = \emptyset$.
Then there is a canonical morphism
\[
cyc: \ch_{m-n}(V) \longto  
\Hom_{\DeffgM} (M (W / Z), \BZ(n)[2n]) \; .
\]
Given an open
immersion $j: U \into W$ and a closed sub-scheme 
$Z'$ of the intersection $Z \cap U$ such that arbitrary intersections of
the irreducible components of $Z'$ are smooth, the diagram
\[
\vcenter{\xymatrix@R-10pt{
        \ch_{m-n}(V) \ar[d]_{j^*} \ar[r]^-{cyc} &
        \Hom_{\DeffgM} (M (W / Z), \BZ(n)[2n]) \ar[d]^-{j^*}  \\
        \ch_{m-n}(V \cap U) \ar[r]^-{cyc} &
        \Hom_{\DeffgM} (M (U / Z'), \BZ(n)[2n]) 
\\}}
\]
commutes.
\end{Cor}

Now fix $X, Y \in PropSm/S$.
Choose a compactification (over $k$)
$\bS$ of $S$, and compactifications 
$\bX$ of $X$, and $\bY$ of $Y$ together with cartesian diagrams
\[
\vcenter{\xymatrix@R-10pt{
        X \ar[d] \ar@{^{ (}->}[r] &
        \bX \ar[d]  \\
        S \ar@{^{ (}->}[r] &
        \bS 
\\}}
\]
and
\[
\vcenter{\xymatrix@R-10pt{
        Y \ar[d] \ar@{^{ (}->}[r] &
        \bY \ar[d]  \\
        S \ar@{^{ (}->}[r] &
        \bS 
\\}}
\]
(this is possible since $X$ and $Y$ are proper over $S$).
The hypothesis on $k$ ensures that 
arbitrary intersections of
the irreducible components of
the complements $\dX$ of $X$ in $\bX$ and
$\dY$ of $Y$ in $\bY$ can be supposed to be smooth.
Each of the three constituents $M$, $\Mcgm$,
$\dMgm$ of the exact triangle $(\ast)$ will correspond to
an application of Corollary~\ref{2bF}, with different
choices of $(W,Z)$.
\begin{enumerate}
\item[(1)] for $M$, we define
$W := X \times_k \bY$,
\item[(2)] for $\Mcgm$, we define
$W := \bX \times_k Y$,
\item[(3)] for $\dMgm$, we define
$W := \bX \times_k \bY - \dX \times_k \dY$.
\end{enumerate}
In all three cases, we put $Z := W - X \times_k Y$. That is,
\begin{enumerate}
\item[(1)] $Z = X \times_k \dY$,
\item[(2)] $Z = \dX \times_k Y$,
\item[(3)] $Z = X \times_k \dY \cup \dX \times_k Y$.
\end{enumerate}
We also let $V := X \times_S Y \subset X \times_k Y$
in all three cases.
These choices satisfy the hypotheses of Corollary~\ref{2bF}
thanks to the following.

\begin{Lem} 
The scheme $X \times_S Y$ is closed in 
$\bX \times_k \bY - \dX \times_k \dY$.
\end{Lem}

\begin{Proof}
Indeed, the diagram
\[
\vcenter{\xymatrix@R-10pt{
        X \times_S Y \ar[d] \ar@{^{ (}->}[r] &
        \bX \times_k \bY - \dX \times_k \dY \ar[d]  \\
        \bS \ar@{^{ (}->}[r]^-{\Delta} &
        \bS \times_k \bS 
\\}}
\]
is cartesian. 
\end{Proof}

\medskip

\begin{Proofof}{Theorem~\ref{2bB}}
We may clearly assume $S$, $X$ and $Y$  
to be of pure absolute dimension $d_S$, $d_X$ and $d_Y$, respectively.

Let us treat $M$ first.
Note that by \cite[Thm.~4.3.7~3]{V}, 
the group of morphisms in $\DgM$ from $M(X)$ to $M(Y)$
is canonically isomorphic to
\[
\Hom_{\DgM} \bigl( M(X) \otimes \Mcgm(Y) , \BZ(d_Y)[2d_Y] \bigr) \; .
\]
Localization for the motive with compact support \cite[Prop.~4.1.5]{V}
shows that $\Mcgm(Y) = M (\bY / \dY)$.
Given the definition of 
the tensor structure on $\DgM$, the above therefore equals
\[
\Hom_{\DgM} \bigl( M(X \times_k \bY / X \times_k \dY), 
\BZ(d_Y)[2d_Y] \bigr) \; .
\]
By Corollary~\ref{2bF}, applied to the setting (1),
this group is the target of the morphism $cyc_1$ on 
$\ch_{d_X}(X \times_S Y) = \ch^{d_Y - d_S}(X \times_S Y)$.
Note that on a class which comes from 
$\FZ \in c_S(X,Y)$, the map $cyc_1$ takes indeed the value $M(\FZ)$.

The case of $\Mcgm$ is similar.
First, by duality,
the group of morphisms in $\DgM$ from $\Mcgm(X)$ to $\Mcgm(Y)$
is canonically isomorphic to
\[
\Hom_{\DgM} \bigl( \Mcgm(X) \otimes M(Y) , \BZ(d_Y)[2d_Y] \bigr) \; .
\]
By localization, this group then equals
\[
\Hom_{\DgM} \bigl( M(\bX \times_k Y / \dX \times_k Y), 
\BZ(d_Y)[2d_Y] \bigr) \; .
\]
By Corollary~\ref{2bF}, applied to the setting (2),
this group is the target of the morphism $cyc_2$ on 
$\ch^{d_Y - d_S}(X \times_S Y)$.

In order to show that for a cycle class $z$ in $\ch^{d_Y - d_S}(X \times_S Y)$,
the diagram
\[
\vcenter{\xymatrix@R-10pt{
        M(X) \ar[d]_{cyc_1(z)} \ar[r] &
        \Mcgm(X) \ar[d]^{cyc_2(z)} \\
        M(Y) \ar[r] &
        \Mcgm(Y) 
\\}}
\]
commutes, we need to study
the group of morphisms in $\DgM$ from $M(X)$ to $\Mcgm(Y)$.
Again by duality, it
is canonically isomorphic to
\[
\Hom_{\DgM} \bigl( M(X \times_k Y), 
\BZ(d_Y)[2d_Y] \bigr) \; .
\]
The above commutativity then follows from the compatibility of $cyc$
under restriction from $X \times_k \bY$, resp.\
$\bX \times_k Y$, to $X \times_k Y$ (Corollary~\ref{2bF}).

Now let us treat $\dMgm$.
Note that by \cite[Thm.~6.1]{W1}, 
the group of morphisms in $\DgM$ from $\dMgm(X)$ to $\dMgm(Y)$
is canonically isomorphic to
\[
\Hom_{\DgM} \bigl( \dMgm(X) \otimes \dMgm(Y)[1] , 
\BZ(d_Y)[2d_Y] \bigr) \; .
\]
As in \cite[pp.~650--651]{W1}, one shows that $\dMgm(X) \otimes \dMgm(Y)[1]$
maps canonically to the relative motive
\[
M \bigl( (\bX \times_k \bY - \dX \times_k \dY) / 
(\bX \times_k \bY - \dX \times_k \dY - X \times_k Y) \bigr) \; .
\]
Hence the group of morphisms 
$\Hom_{\DeffQgM} ( \dMgm(X) , \dMgm(Y))$
receives an arrow, say $\alpha$, from the group of morphisms from
\[
M \bigl( (\bX \times_k \bY - \dX \times_k \dY) / 
(\bX \times_k \bY - \dX \times_k \dY - X \times_k Y) \bigr)
\]
to $\BZ(d_Y)[2d_Y])$.
By Corollary~\ref{2bF}, applied to the setting (3),
this group is the target of the morphism $cyc_3$ on 
$\ch^{d_Y - d_S}(X \times_S Y)$.

In order to show that for a cycle class $z$ in $\ch^{d_Y - d_S}(X \times_S Y)$,
the diagram
\[
\vcenter{\xymatrix@R-10pt{
        \Mcgm(X) \ar[d]_{cyc_2(z)} \ar[r] &
        \dMgm(X)[1] \ar[d]^{cyc_3(z)[1]} \\
        \Mcgm(Y) \ar[r] &
        \dMgm(Y)[1] 
\\}}
\]
commutes, we need to study
the group of morphisms in $\DgM$ from $\Mcgm(X)$ to $\dMgm(Y)[1]$.
Again by \cite[Thm.~6.1]{W1}, it
is canonically isomorphic to
\[
\Hom_{\DgM} \bigl( \Mcgm(X) \otimes \dMgm(Y), 
\BZ(d_Y)[2d_Y] \bigr) \; .
\]
But $\Mcgm(X) \otimes \dMgm(Y)$
maps canonically to $\Mcgm(X) \otimes M(Y)$, 
which was already identified with the relative motive
\[
M(\bX \times_k Y / \dX \times_k Y) \; .
\]
Hence the group of morphisms 
$\Hom_{\DgM} (\Mcgm(X) , \dMgm(Y)[1])$ receives an arrow, say $\beta$, from
\[
\Hom_{\DgM} \bigl( M(\bX \times_k Y / \dX \times_k Y), 
\BZ(d_Y)[2d_Y] \bigr) \; .
\]
The desired commutativity then follows from the compatibility of $cyc$
under restriction from $\bX \times_k \bY - \dX \times_k \dY$ 
to $\bX \times_k Y$ (Corollary~\ref{2bF}),
and from the compatibility of $\beta$ and the map $\alpha$ from above.
The latter is a consequence of the compatibility of the isomorphism
\[
\dMgm(Y)[1] \isoto \dMgm(Y)^*(d_Y)[2d_Y]
\]
with duality $M(Y) \cong \Mcgm(Y)^*(d_Y)[2d_Y]$
\cite[Thm.~6.1]{W1}.

The proof of the commutativity of
\[
\vcenter{\xymatrix@R-10pt{
        \dMgm(X) \ar[d]_{cyc_3(z)} \ar[r] &
        M(X) \ar[d]^{cyc_1(z)} \\
        \dMgm(Y) \ar[r] &
        M(Y) 
\\}}
\]
is similar.

Altogether, this proves part~(a) of the statement.
As for parts~(b) and (c), simply compose the functor from (a) with
duality in $\DgM$, using \cite[Thm.~4.3.7~3]{V} and \cite[Thm.~6.1]{W1}.

By \cite[Rem.~6.15]{W1}, our construction
is independent of the compactifications $\bS$, $\bX$, $\bY$.
\end{Proofof}

\medskip

\begin{Proofof}{Theorem~\ref{2bBa}}
We keep the notations of the previous proof. Choose
compactifications 
$\bT$ of $T$, and $\bU$ of $U$ together with cartesian diagrams
\[
\vcenter{\xymatrix@R-10pt{
        T \ar[d]_f \ar@{^{ (}->}[r] &
        \bT \ar[d]  \\
        S \ar@{^{ (}->}[r] &
        \bS 
\\}}
\]
and
\[
\vcenter{\xymatrix@R-10pt{
        U \ar[d]_g \ar@{^{ (}->}[r] &
        \bU \ar[d]  \\
        S \ar@{^{ (}->}[r] &
        \bS 
\\}}
\]
($f$ and $g$ are proper).

(a)~Checking the definitions, the transformation $\alpha_{f_\sharp}$
is in fact given by the identity. Indeed, both
$( \dMgm , M , \Mcgm )_S \circ f_\sharp$ and $( \dMgm , M , \Mcgm )_T$ 
map the object $h(X/T)$, for $X \in PropSm/T$, to the exact triangle
\[
(\ast)_X \quad\quad
\dMgm(X) \longto M(X) \longto \Mcgm(X) \longto \dMgm(X)[1] \; .
\]
Note that on morphisms, the functor $f_\sharp$ corresponds
to the push-forward 
\[
\ch_*(X \times_T Y) \longto \ch_*(X \times_S Y)
\]
along the closed immersion $X \times_T Y \into X \times_S Y$.
The latter factors the closed immersion 
\[
X \times_T Y \longinto \bX \times_k \bY - \dX \times_k \dY \; .
\]
The construction (see the preceding proof) shows then that
the effects of the functors $( \dMgm , M , \Mcgm )_S \circ f_\sharp$ and
of $( \dMgm , M , \Mcgm )_T$ coincide also
on $\ch_*(X \times_T Y)$.
This shows the first half of statement~(a). The second is implied
formally by Theorem~\ref{2bB}~(c).
 
(b)~We first consider an auxiliary functor.
The morphism $g$ being proper and smooth, we may consider the
composition $g_\sharp \circ g^*$ on $\CHSM$, which on objects is given by
$h(Y/S) \mapsto h(Y \times_S U/S)$, for any proper smooth
scheme $Y$ over $S$. Projection onto the first component then yields
a transformation of functors, namely the adjunction
\[
b_{g^*,\id_S} : g_\sharp \circ g^* \longto \id_{\CHSM} \; .
\]
Then define $\beta_{g^*,\id_S}$ to be the composition of transformations
\[
\beta_{g^*,\id_S} := 
\bigl( ( \dMgm , M , \Mcgm )_S \circ b_{g^*,\id_S} \bigr) \circ
\bigl( \alpha_{g_{\sharp}} \circ g^* \bigr)^{-1} \; , 
\]
observing the equivalence 
\[
\alpha_{g_{\sharp}} \circ g^*: 
( \dMgm , M , \Mcgm )_S \circ g_\sharp \circ g^* \isoto 
( \dMgm , M , \Mcgm )_U \circ g^* 
\]
from part~(a). 
We leave it to the reader to check the compatibility
of this construction with composition of proper smooth morphisms in $Sm/k$.

(c)~Similarly, this commutativity statement is left as an exercice.

(d), (e)~Given Theorem~\ref{2bB}~(c), these statements follow formally from
(b) and (c), respectively.
\end{Proofof}

For $X,Y \in PropSm/S$, denote by $\bar{c}_S(X,Y)$ the quotient
of $c_S(X,Y)$ by the group of cycles $\FZ$ satisfying
\[
M(\FZ) = 0 \quad \; , \quad 
\Mcgm(\FZ) = 0 \quad \; , \quad 
\dMgm(\FZ) = 0 \; .
\]
Note that composition of correspondences induces a well-defined
composition on $\bar{c}_S$. In particular, for any $X \in PropSm/S$,
the group $\bar{c}_S(X,X)$ carries the structure of an algebra.

\begin{Cor} \label{2bG} 
Let $X$ and $Y$ be in $PropSm/S$. 
Then the projection
\[
c_S(X,Y) \longonto \bar{c}_S(X,Y)
\]
factors through the image of $c_S(X,Y)$ in 
$\ch^* (X \times_S Y)$. 
In other words, two cycles $\FZ_1, \FZ_2 \in c_S(X,Y)$
induce the same morphisms $M(\FZ_i)$, resp.\ $\Mcgm(\FZ_i)$,
resp.\ $\dMgm(\FZ_i)$,
if they are rationally equivalent (on $X \times_S Y$).
\end{Cor}

\forget{
\begin{Cor} \label{2bH}
Let $F$ be a flat $\BZ$-algebra,
$X \in PropSm/S$, $I$ a set of indices 
and $\FZ_i$, $i \in I$, cycles in $c_S(X,X) \otimes_\BZ F$.
Denote by $\bar{\FZ}_i$ their classes in $\bar{c}_S(X,X) \otimes_\BZ F$. 
\\[0.1cm] 
(a)~If the images of the $\FZ_i$ in 
$\ch^* (X \times_S X) \otimes_\BZ F$ are idempotent,
then so are the $\bar{\FZ}_i$. \\[0.1cm]
(b)~If the images of the $\FZ_i$ in 
$\ch^* (X \times_S X) \otimes_\BZ F$ commute with each other,
then so do the $\bar{\FZ}_i$.
\end{Cor}
}
\begin{Rem} \label{2bGa}
(a)~Let $X,Y \in Sm/S$.
As shown in \cite[Lemma~5.18]{L3}, the map
\[
c_S(X,Y) \longto \ch_{d_X} (X \times_S Y)
\]
is surjective, whenever $Y$ is projective, and $X$
of pure absolute dimension $d_X$. Therefore, 
by Corollary~\ref{2bG}, the group
$\bar{c}_S(X,Y)$ is canonically a quotient of $\ch_{d_X} (X \times_S Y)$
if $X \in PropSm/S$ and $Y \in ProjSm/S$. \\[0.1cm]
(b)~The observation from (a)
fits in the functorial picture sketched in Remark~\ref{2bC}.
Indeed, \cite[Thm.~5.23]{Deg2} implies that 
the restriction of the functor $M$,
\[
M: \CHSM_{proj} \longto \DgM
\]
factors canonically through a fully faithful embedding 
\[
\CHSM_{proj} \longinto \DSgM \; .
\]
(c)~In \cite[Prop.~5.19]{L3}, an  
embedding result analogous to~(b) is proven for a dg-version
of $\DSgM$, denoted $SmMot(S)$ in \loccit , from which 
the embedding~(b)
can be deduced \cite[Cor.~7.13]{L3}. \\[0.1cm]
(d)~Using \cite[Thm.~5.23]{Deg2}, one can show that 
\[
\Hom_{\DSgM} ( M_1 ,M_2[i] ) = 0
\]
for any two smooth relative Chow motives $M_1, M_2 \in \CHSM_{proj}$, 
and any integer $i > 0$. \\[0.1cm]
(e)~When $S = \Spec k$, results~(b) and (d) are 
contained in \cite[Cor.~4.2.6]{V}.
\end{Rem}

\begin{Rem} \label{2bI}
Fix a non-negative integer $d$, and consider the full sub-category
$\CHSM_d$ of $\CHSM$ of smooth relative Chow motives generated by the Tate twists
of $h(X/S)$, for $X \in PropSm/S$ of pure absolute dimension $d$.
The construction of the duality isomorphisms 
\cite[Thm.~4.3.7~3]{V}, \cite[Thm.~6.1]{W1} shows that the identification
\[
\bigl( \dMgm , M , \Mcgm \bigr)^* = 
\bigl( \dMgm , M , \Mcgm \bigr)(-d)[-2d]
\]
of the restriction of the functors from Theorem~\ref{2bB} to $\CHSM_d$
admits an alternative description, when $S$ is of pure absolute dimension,
say $s$: on $\CHSM_d$, the functor $( \dMgm , M , \Mcgm )^*$
equals then composition of duality in the cate\-gory $\CHSM_d$ with 
$( \dMgm , M , \Mcgm )$, followed by the functor $M \mapsto M(-s)[-2s]$.
Note that on morphisms, duality in $\CHSM_d$ corresponds to
the transposition $\ch_* (X \times_S Y) \to \ch_* (Y \times_S X)$.
\end{Rem}

This observation allows to deduce the following statements from
Theorem~\ref{2bBa}.

\begin{Cor} \label{2bK}
(a)~Let $g: U \to S$ be a proper smooth morphism in $Sm/k$
of pure relative dimension $d_g$.
Then there exists a canonical transformation of additive functors
\[
\delta_{\id_S,g^*} : ( \dMgm , M , \Mcgm )_S \longto 
( \dMgm , M , \Mcgm )_U (-d_g)[-2d_g] \circ g^* \; .
\]
The formation of $\delta_{\id_S,g^*}$ is compatible with composition
of proper smooth morphisms in $Sm/k$ of pure relative dimension.  \\[0.1cm]
(b)~Let $g: U \to S$ be a finite \'etale morphism in $Sm/k$ 
of constant (fibrewise) degree $u$. 
Then the endomorphism 
\[
\beta_{g^*,\id_S} \circ \delta_{\id_S,g^*} 
\]
of the functor $( \dMgm , M , \Mcgm )_S$ equals
multiplication by $u$.
\end{Cor}

\begin{Proof}
(a)~We may assume $S$ to be of pure absolute dimension, say $s$.
Consider the transformation 
\[
\gamma_{\id_S,g^*} : 
( \dMgm , M , \Mcgm )^*_S \longto 
( \dMgm , M , \Mcgm )^*_U \circ g^* 
\]
from Theorem~\ref{2bBa}~(d). Composition with duality $\BD_S$ in $\CHSM$ gives
\[
\gamma_{\id_S,g^*} \circ \BD_S: 
( \dMgm , M , \Mcgm )^*_S \circ \BD_S \longto 
( \dMgm , M , \Mcgm )^*_U \circ g^* \circ \BD_S \; . 
\]
Now observe the formula 
\[
\BD_U \circ g^* = g^* \circ \BD_S
\]
($\BD_U :=$ duality in $\CHUM$). Define 
$\delta_{\id_S,g^*}$ as the composition of $\gamma_{\id_S,g^*} \circ \BD_S$
and $M \mapsto M(s)[2s]$, observing that
source and target of $\delta_{\id_S,g^*}$ are identified with 
$( \dMgm , M , \Mcgm )_S$ and  
$( \dMgm , M , \Mcgm )_U (-d_g)[-2d_g] \circ g^*$,
respectively.

(b)~The morphism $g$ being finite and \'etale, we have
\[
\BD_S \circ g_\sharp = g_\sharp \circ \BD_U \; .
\]
This shows that $g_\sharp$ is also right adjoint to $g^*$.
Checking the definitions, the composition 
$\beta_{g^*,\id_S} \circ \delta_{\id_S,g^*}$ equals up to twist and shift
the composition of the two adjunctions
\[
\xi: \id_{\CHSM} \longto g_\sharp \circ g^* \longto \id_{\CHSM} \; ,
\]
preceded by duality, and followed by $( \dMgm , M , \Mcgm )^*_S$.
These functors being additive, it suffices to show that $\xi$
equals multiplication bu $u$. But this identity 
on morphisms of smooth relative Chow motives is classical.
\end{Proof}

The main results of this section have obvious $F$-linear versions,
for any commutative $\BQ$-algebra $F$.
Let us now describe how our analysis of the 
functor $( \dMgm , M , \Mcgm )$
will be used in the sequel.

\begin{Ex} \label{2bL}
Let $g_1, g_2: U \to S$ be two finite \'etale morphisms in $Sm/k$.
Fix an object $X \in PropSm/S$, an idempotent $e$ on
$h(X/S)$ (possibly belonging to $\ch^*(X \times_S X) \otimes_\BZ F$,
for some commutative $\BQ$-algebra $F$), and a morphism 
\[
\varphi: g_1^* \bigl( h(X/S)^e \bigr) \longto g_2^* \bigl( h(X/S)^e \bigr) 
\]
in $\CHUM$ (or $\CHUM_F$). \\[0.1cm]
(a)~Let us define an endomorphism of
$( \dMgm , M , \Mcgm )(h(X/S)^e)$ ``of Hecke type'',
denoted $\varphi(g_1,g_2)$, by composing
\[
\delta_{\id_S,g_1^*} : (\dMgm , M , \Mcgm) \bigl( h(X/S)^e \bigr) \longto 
( \dMgm , M , \Mcgm ) \bigl( g_1^* \bigl( h(X/S)^e \bigr) \bigr) 
\]
first with $( \dMgm , M , \Mcgm ) \circ \varphi$, and then with
\[
\beta_{g_2^*,\id_S} : 
( \dMgm , M , \Mcgm )  \bigl( g_2^* \bigl( h(X/S)^e \bigr) \bigr) \to
( \dMgm , M , \Mcgm ) \bigl( h(X/S)^e \bigr) \; .
\]
(b)~Note that unless $g_1 = g_2$, the endomorphism 
$\varphi(g_1,g_2)$ is in general \emph{not}
the image of an endomorphism on the smooth relative Chow motive $h(X/S)^e$
under the functor $( \dMgm , M , \Mcgm )$. \\[0.1cm]
(c)~If $\varphi$ is an isomorphism, 
with inverse $\psi$, then using the construction
from (a), the endomorphism 
$\psi(g_2,g_1)$ on
$( \dMgm , M , \Mcgm )(h(X/S)^e)$ can be defined. If $X$
is of pure absolute dimension $d_X$, then 
$\psi(g_2,g_1)$ equals the dual of $\varphi(g_1,g_2)$, 
twisted by $d_X$ and shifted by $2d_X$, under the identification  
\[
( \dMgm , M , \Mcgm )^* \bigl( h(X/S) \bigr) =
( \dMgm , M , \Mcgm ) \bigl( h(X/S) \bigr)(-d_X)[-2d_X] 
\]
from Theorem~\ref{2bB}~(b1). We leave the details of the verification
to the reader. \\[0.1cm]
(d)~In practice, the morphism 
$\varphi: g_1^* ( h(X/S)^e ) \to g_2^* ( h(X/S)^e )$ 
will be obtained
from a morphism of smooth relative Chow motives over $U$
\[
\varphi: h(X \times_{S,g_1} U / U) = g_1^* ( h(X/S) ) \longto 
g_2^* ( h(X/S) ) = h(X \times_{S,g_2} U / U) 
\]
satisfying the equation 
\[
\varphi \circ g_1^*(e) = g_2^*(e) \circ \varphi
\]
in $\ch^* ( (X \times_{S,g_1} U) \times_U (X \times_{S,g_2} U) )$ 
(or $\ch^* ( (X \times_{S,g_1} U) \times_U (X \times_{S,g_2} U) )
\otimes_\BZ F$). In that case, $\varphi(g_1,g_2)$ can be seen as an
endomorphism of the whole of $( \dMgm , M , \Mcgm )(h(X/S))$
commuting with $e$. \\[0.1cm]
(e)~In the setting of (d), assume that the morphism
\[
\varphi: h(X \times_{S,g_1} U / U) \longto h(X \times_{S,g_2} U / U) 
\]
is represented by the cycle $\FZ$ in
\[
c_U (X \times_{S,g_1} U , X \times_{S,g_2} U) 
\]
(or in $c_U (X \times_{S,g_1} U , X \times_{S,g_2} U) \otimes_\BZ F$).
Checking the definitions, one sees that
the $M$-component of $\varphi(g_1,g_2)$ is then represented
by the image of $\FZ$ under the direct image 
\[
\bigl( g_1 \times_k g_2 \bigr)_* : 
c_U (X \times_{S,g_1} U , X \times_{S,g_2} U) \longto c(X,X) 
\]
(or under $( g_1 \times_k g_2 )_* \otimes F$).
\end{Ex}

%%% Local Variables:
%%% mode: latex
%%% TeX-master: "head"
%%% End:

\bigskip
%\include{Sec2c}

%%%%%%%%%%%%%%%%%%%%%%%%%%%%%%%%%%%%%%%%%%%%%%%%%%%%%%%%%%%%%%%%%%%%%%%
%
%  Section 2c
%
%%%%%%%%%%%%%%%%%%%%%%%%%%%%%%%%%%%%%%%%%%%%%%%%%%%%%%%%%%%%%%%%%%%%%%%

\section{Motives associated to Abelian schemes}
\label{2c}

%%%%%%%%%%%%%%%%%%%%%%%%%%%%%%%%

%%%%%%%%%%%%%%%%%%%%%%%%%%%%%%%%

Fix a field $k$ admitting
strict resolution of singularities, and a base $S \in Sm/S$. 
In this section, we combine the main result
from \cite{DM} with the theory developed in Section~\ref{2b}. 
Recall the following.

\begin{Thm}[{\cite[Thm.~3.1, Prop.~3.3]{DM}}] \label{2cA}
(a)~Let $A/S$ be an \\ Ab\-elian scheme of relative dimension $g$.
Then there is a unique decomposition of the class of the diagonal 
$(\Delta) \in \ch^g (A \times_S A) \otimes_\BZ \BQ \;$,
\[
(\Delta) = \sum_{i=0}^{2g} p_{A,i}
\]
such that $p_{A,i} \circ (\Gamma_{[n]_A}) = n^i \cdot p_{A,i}$ for all $i$,
and all integers $n$.  
The $p_{A,i}$ are mutually orthogonal idempotents, and
$(\Gamma_{[n]_A}) \circ p_{A,i} = n^i \cdot p_{A,i}$ for all $i$. \\[0.1cm]
(b)~For any morphism $f: A \to B$ of Abelian schemes over $S$, and any $i$, 
\[
p_{B,i} \circ (\Gamma_{f}) = (\Gamma_{f}) \circ p_{A,i}
\in \ch^* (A \times_S B) \otimes_\BZ \BQ \; .
\]
In other words, the decomposition in (a) is 
covariantly functorial in $A$. \\[0.1cm]
(c)~For any isogeny $g: B \to A$ of Abelian schemes over $S$, and any $i$, 
\[
p_{B,i} \circ ({}^t \Gamma_{g}) = ({}^t \Gamma_{g}) \circ p_{A,i}
\in \ch^* (A \times_S B) \otimes_\BZ \BQ \; .
\]
In other words, the decomposition in (a) is contravariantly 
functorial under isogenies.
\end{Thm}

We use the notation $\Gamma_h$ for the graph of a morphism
$h$ of $S$-schemes, $[n]_A$ for the multiplication by $n$
on the Abelian scheme $A$, $(\FZ)$ for the class of a cycle $\FZ$,
and ${}^t \FZ$ for its transposition. 
Let 
\[
h(A/S) = \bigoplus_i h_i(A/S)
\]
be the decomposition of the relative motive of $A$
corresponding to the decomposition $(\Delta) = \sum_i p_{A,i}$.
Thus, on the term $h_i(A/S)$, the cycle class $(\Gamma_{[n]_A})$
acts \emph{via} multiplication by $n^i$. \\

Now recall the exact triangle
\[
(\ast)_A \quad\quad
\dMgm(A) \longto M(A) \longto \Mcgm(A) \longto \dMgm(A)[1] \; .
\]
By Theorem~\ref{2bB}~(a), the cycle classes $p_{A,i}$
induce endomorphisms of $(\ast)_A$, 
when considered as an exact triangle in $\DeffqgM$.

\begin{Thm} \label{2cB} 
(a)~Let $A/S$ be an Abelian scheme
of relative dimension $g$.
For $0 \le i \le 2g$, denote by $M(A)_i$,
$\Mcgm(A)_i$ and $\dMgm(A)_i$ the images of the idempotent
$p_{A,i}$ on $M(A)$, $\Mcgm(A)$ and $\dMgm(A)$,
respectively, considered as objects of the category $\DeffqgM$.
Then for any $i$, the triangle
\[
(\ast)_{A,i} \quad\quad
\dMgm(A)_i \longto M(A)_i \longto \Mcgm(A)_i \longto \dMgm(A)_i[1] 
\]
in $\DeffqgM$ is exact. \\[0.1cm]
(b)~The direct sum of the triangles $(\ast)_{A,i}$ 
yields a decomposition
\[
(\ast)_A = \; \bigoplus_{i=0}^{2g} \; \; (\ast)_{A,i} \; .
\]
It has the following properties:
\begin{enumerate}
\item[(b1)] for any integer $n$, the decomposition is respected
by $[n]_A$. 
\item[(b2)] for each $i$ and $n$, the induced morphisms $[n]_{A,i}$ on the 
three terms of $(\ast)_{A,i}$ equal multiplication by $n^i$.  
\end{enumerate}
(c)~As a decomposition of $(\ast)_A$ into 
some finite direct sum of exact triangles
in $\DeffqgM$, 
\[
(\ast)_A = \; \bigoplus_i \; \; (\ast)_{A,i}
\]
is uniquely determined by properties~(b1) and (b2). More precisely,
it is uniquely determined by the following properties:
\begin{enumerate}
\item[(c1)] for some integer $n \ne -1,0,1$, 
the decomposition is respected
by $[n]_A$. 
\item[(c2)] for the choice of $n$ made in~(c1) and each $i$, 
the induced morphism $[n]_{A,i}$ on the 
three terms of $(\ast)_{A,i}$ equals multiplication by $n^i$.  
\end{enumerate}
(d)~The decomposition 
\[
(\ast)_A = \; \bigoplus_i \; \; (\ast)_{A,i}
\]
is covariantly functorial under morphisms,
and contravariantly functorial under isogenies of Abelian schemes over $S$.
\end{Thm}

\begin{Proof}
Part~(a) is a formal consequence of the fact that the $p_{A,i}$
are idempotent.

Parts~(b) and (d) follow from Theorem~\ref{2cA}
and the functoriality statement from Theorem~\ref{2bB}~(a).

Part~(c) is left to the reader.
\end{Proof}

The following seems worthwhile to note explicitly.

\begin{Cor} \label{2cC}
Let $A/S$ be an Abelian scheme of relative dimension $g$.
Then the boundary motive $\dMgm(A)$ decomposes functorially into a direct sum
\[
\dMgm(A) = \bigoplus_{i=0}^{2g} \dMgm(A)_i \; .
\]
On $\dMgm(A)_i$, the endomorphism $[n]_A$ 
acts \emph{via} multiplication by $n^i$, for any integer $n$,
and any $0 \le i \le 2g$.
\end{Cor}

Here is an illustration
of the surjectivity proved in \cite[Lemma~5.18]{L3}. 
 
\begin{Prop} \label{2cD}
Let $A/S$ be an Ab\-elian scheme.
The elements $p_{A,i}$ of $\ch^* (A \times_S A) \otimes_\BZ \BQ$
lie in the image of 
\[
c_S(A,A) \otimes_\BZ \BQ \longto 
\ch^* (A \times_S A) \otimes_\BZ \BQ \; .
\]
More precisely, for any integer $n \ne -1, 0, 1$,
\[
\pi_{A,i,n} := \prod_{j \ne i} \frac{\Gamma_{[n]_A} - n^j}{n^i - n^j}
\]
is a pre-image of $p_{A,i}$ in $c_S(A,A) \otimes_\BZ \BQ \;$. 
\end{Prop}

\begin{Proof}
On each of the direct factors $h_j(A/S) \subset h(A/S)$, 
the projector $p_{A,i}$ acts \emph{via} 
multiplication by the Kronecker symbol $\delta_{ij}$, while $(\Gamma_{[n]_A})$
acts \emph{via} multiplication by $n^j$. Therefore,
\[
p_{A,i} = \prod_{j \ne i} \frac{(\Gamma_{[n]_A}) - n^j}{n^i - n^j} 
\in \ch^* (A \times_S A) \otimes_\BZ \BQ \; ,
\]
for any integer $n \ne -1, 0, 1$.
Therefore, the element
$\pi_{A,i,n} \in c_S(A,A) \otimes_\BZ \BQ$
is indeed a pre-image of $p_{A,i}$.  
\end{Proof}

\bigskip

%\include{Sec2}

%%%%%%%%%%%%%%%%%%%%%%%%%%%%%%%%%%%%%%%%%%%%%%%%%%%%%%%%%%%%%%%%%%%%%%%
%
%  Section 2
%
%%%%%%%%%%%%%%%%%%%%%%%%%%%%%%%%%%%%%%%%%%%%%%%%%%%%%%%%%%%%%%%%%%%%%%%

\section{The intersection motive of a surface}
\label{2}

%%%%%%%%%%%%%%%%%%%%%%%%%%%%%%%%

%%%%%%%%%%%%%%%%%%%%%%%%%%%%%%%%

Fix a normal, proper surface $X^*$ over $k$. 
Let us first recall, following \cite{CM}, the construction 
and the basic properties of the 
\emph{intersection motive} of $X^*$. Choose
\[
\vcenter{\xymatrix@R-10pt{
        X \ar@{^{ (}->}[r] &
        X^* \ar@{<-^{ )}}[r] &
        Z
\\}}
\]
where $Z$ is a closed sub-scheme of $X^*$ which is finite over $k$,
and whose complement $X$ is  
smooth. Choose a
resolution of singularities. More precisely, consider in addition 
the following diagram, assumed
to be cartesian:
\[
\vcenter{\xymatrix@R-10pt{
        X \ar@{^{ (}->}[r] \ar@{=}[d] &
        {\Xp} \ar@{<-^{ )}}[r]^{i} \ar[d]_\pi &
        D \ar[d]^\pi \\
        X \ar@{^{ (}->}[r] &
        X^* \ar@{<-^{ )}}[r] &
        Z
\\}}
\]
where $\pi$ is proper (and birational), $\Xp$ is smooth (and proper), 
and $D$ is a divisor with
normal crossings, whose irreducible components $D_m$ are smooth (and proper). \\

Recall \cite[Sect.~1.13]{Sch} that the ``degree $2$ parts''
$M_2(D_m)$ are canonically defined as sub-objects of the motives 
$M(D_m)$ (we remind the reader that throughout the article, we
use homological notation). Hence there is a canonical morphism
\[
i_{*,2}: M_2(D) := \bigoplus_m M_2(D_m) \longinto \bigoplus_m M(D_m) 
\longto M(\Xp)  
\]
of Chow motives. Similarly \cite[Sect.~1.11]{Sch}, 
there is a canonical morphism
\[
i^*_0: M(\Xp) \longto \bigoplus_m M(D_m)(1)[2] 
\longonto \bigoplus_m M_0(D_m)(1)[2]   \; ,
\]
where $M_0(D_m)$ denote the ``degree $0$ parts'',
canonically defined as quotients of $M(D_m)$. 
The following is a special case of \cite[Sect.~2.5]{CM}
(see also \cite[Thm.~2.2]{W4}).

\begin{Thm} \label{2A}
(i) The composition $\alpha := i^*_0i_{*,2}$ is an isomorphism
in the $\BQ$-linear category $\DeffqgM$. \\[0.1cm]
(ii) The composition $p:= i_{*,2}\alpha^{-1}i^*_0$ is an idempotent on 
$M(\Xp) \in \DeffqgM$. 
Hence so is the difference $\id_{\Xp}-p$. \\[0.1cm]
(iii) The image $\imm p \in \DeffqgM$ 
is canonically isomorphic to $M_2(D)$.
\end{Thm}

The proof relies on the non-degeneracy of the intersection pairing
on the components of $D$. 

\begin{Def}[{\cite[p.~158]{CM}, \cite[Def.~2.3]{W4}}] \label{2B}
The intersection \\
motive of $X^*$ is defined as
\[
M^{!*} (X^*) := \imm (\id_{\Xp}-p) \in \DeffqgM \; .
\]
\end{Def}

The name is motivated by the behaviour of the realizations
of the intersection motive. Its
functoriality properties are
given in \cite[Prop.~2.5]{W4}. It will be useful to recall 
in particular the behaviour
under finite morphisms
$f: Y^* \to X^*$ between
normal, proper surfaces over $k$.
Assume that $Z$ is such that the
pre-image under $f$ of $X = X^* - Z$ is dense, and smooth
(this can be achieved by enlarging $Z$, if necessary). 
The closed sub-scheme $f^{-1}(Z)$ of
$Y$ contains the singularities of $Y^*$. 
We thus can find a cartesian diagram of desingularizations of
$X^*$ and $Y^*$ of the type considered before:
\[
\vcenter{\xymatrix@R-10pt{
        \Yp \ar@{<-^{ )}}[r] \ar[d]_F &
        C \ar[d]^F \\
        \Xp \ar@{<-^{ )}}[r] &
        D
\\}}
\]
The following is the content of \cite[Prop.~2.5~(iii) and (iv), Prop.~2.4]{W4}.

\begin{Prop} \label{2C}
(a)~Both $F^*$ and $F_*$ respect the decompositions
\[
M(\Xp) = M^{!*}(X^*) \oplus M_2(D)
\]
and
\[
M(\Yp) = M^{!*}(Y^*) \oplus M_2(C)
\]
of $M(\Xp)$ and of $M(\Yp)$, respectively.
The composition $F_*F^*$ equals multiplication with the degree of $f$. \\[0.1cm]
(b)~The definition of $M^{!*} (X^*)$ is
independent of the choices of the finite sub-scheme 
$Z$ containing the singularities, and of
the desingularization $\Xp$ of $X^*$.
\end{Prop}

Next, let us establish the connection to the boundary motive of $X$,
and to the constructions of Section~\ref{1}. To do so, 
assume $k$ to admit resolution of singularities, \emph{fix}
a dense open sub-scheme $X \subset X^*$ which is smooth, and \emph{choose}
\[
\vcenter{\xymatrix@R-10pt{
        X \ar@{^{ (}->}[r] \ar@{=}[d] &
        {\Xp} \ar@{<-^{ )}}[r]^{i} \ar[d]_\pi &
        D \ar[d]^\pi \\
        X \ar@{^{ (}->}[r] &
        X^* \ar@{<-^{ )}}[r] &
        Z
\\}}
\]
as above. Recall the diagram of exact triangles
\[
\vcenter{\xymatrix@R-10pt{
        0 &
        M(D)^*(2)[4] \ar[l] &
        M(D)^*(2)[4] \ar@{=}[l]  &
        0 \ar[l] \\
        \Mcgm(X) \ar[u] &
        M(\Xp) \ar[u]^{i^*} \ar[l] &
        M(D) \ar[u] \ar[l]_{i_*} &
        \Mcgm(X)[-1] \ar[u] \ar[l] \\
        \Mcgm(X) \ar@{=}[u]&  
        M(X) \ar[u] \ar[l] &
        \dMgm(X) \ar[u] \ar[l] &
        \Mcgm(X)[-1] \ar@{=}[u] \ar[l] \\
        0 \ar[u] &
        M(D)^*(2)[3] \ar[u] \ar[l] &
        M(D)^*(2)[3] \ar[u] \ar@{=}[l] &
        0 \ar[l] \ar[u]   
\\}}
\]
from Theorem~\ref{1N}~(c), and let us refer to it using the symbol
$(A)$. 
It turns out that the three components of the idempotent
$p= i_{*,2}\alpha^{-1}i^*_0$ on $M(\Xp)$ all extend to 
give morphisms of diagrams of exact triangles: the first,
denoted $i^*_0$, maps $(A)$ to  
\[
\vcenter{\xymatrix@R-10pt{
        0 &
        \bigoplus_m M_0(D_m)(1)[2] \ar[l] &
        \bigoplus_m M_0(D_m)(1)[2] \ar@{=}[l]  &
        0 \ar[l] \\
        0 \ar@{=}[u] &
        \bigoplus_m M_0(D_m)(1)[2] \ar@{=}[u] \ar[l] &
        \bigoplus_m M_0(D_m)(1)[2] \ar@{=}[u] \ar@{=}[l] &
        0 \ar@{=}[u] \ar[l] \\
        0 \ar@{=}[u]&  
        0 \ar[u] \ar@{=}[l] &
        0 \ar[u] \ar@{=}[l] &
        0 \ar@{=}[u] \ar@{=}[l] \\
        0 \ar@{=}[u] &
        \bigoplus_m M_0(D_m)(1)[1] \ar[u] \ar[l] &
        \bigoplus_m M_0(D_m)(1)[1] \ar[u] \ar@{=}[l] &
        0 \ar[l] \ar@{=}[u]   
\\}}
\]
The second component $\alpha^{-1}$
maps this diagram isomorphically to the following, which we shall
denote by $(B)$.
\[
\vcenter{\xymatrix@R-10pt{
        0 &
        M_2(D) \ar[l] &
        M_2(D) \ar@{=}[l]  &
        0 \ar[l] \\
        0 \ar@{=}[u] &
        M_2(D) \ar@{=}[u] \ar[l] &
        M_2(D) \ar@{=}[u] \ar@{=}[l] &
        0 \ar@{=}[u] \ar[l] \\
        0 \ar@{=}[u]&  
        0 \ar[u] \ar@{=}[l] &
        0 \ar[u] \ar@{=}[l] &
        0 \ar@{=}[u] \ar@{=}[l] \\
        0 \ar@{=}[u] &
        M_2(D)[-1] \ar[u] \ar[l] &
        M_2(D)[-1] \ar[u] \ar@{=}[l] &
        0 \ar[l] \ar@{=}[u]   
\\}}
\]
Finally the third component $i_{*,2}$ maps $(B)$ back to $(A)$.
The composition of the three morphisms, denoted by
\[
p : (A) \longto (A) \; ,
\]
is idempotent. Its image is diagram $(B)$.
Denote the image of $\id - p$ on $M(D)$ by $M_{\le 1}(D)$.
Then the image of $\id - p$ on the whole diagram equals
\[
\vcenter{\xymatrix@R-10pt{
        0 &
        M_{\le 1}(D)^*(2)[4] \ar[l] &
        M_{\le 1}(D)^*(2)[4] \ar@{=}[l]  &
        0 \ar[l] \\
        \Mcgm(X) \ar[u] &
        M^{!*}(X^*) \ar[u]^{i^*} \ar[l] &
        M_{\le 1}(D) \ar[u] \ar[l]_{i_*} &
        \Mcgm(X)[-1] \ar[u] \ar[l] \\
        \Mcgm(X) \ar@{=}[u]&  
        M(X) \ar[u] \ar[l] &
        \dMgm(X) \ar[u] \ar[l] &
        \Mcgm(X)[-1] \ar@{=}[u] \ar[l] \\
        0 \ar[u] &
        M_{\le 1}(D)^*(2)[3] \ar[u] \ar[l] &
        M_{\le 1}(D)^*(2)[3] \ar[u] \ar@{=}[l] &
        0 \ar[l] \ar[u]   
\\}}
\]

\begin{Thm} \label{2D}
Assume $k$ to admit resolution of singularities. \\[0.1cm]
(a)~For fixed $X^*$ and smooth $X \subset X^*$, the above diagram 
of exact triangles is independent of the choice of $\Xp$. \\[0.1cm]
(b)~The diagram is covariantly and contravariantly functorial 
under finite morphisms $f: Y^* \to X^*$ of normal, proper surfaces,
which are compatible with the choices of smooth sub-schemes
$X \subset X^*$ and $Y \subset Y^*$: $f^{-1}(X) = Y$. \\[0.1cm]
(c)~The third column of the diagram
\[
M_{\le 1}(D)^*(2)[3] \longto \dMgm(X) \longto M_{\le 1}(D)
\longto M_{\le 1}(D)^*(2)[4]
\]
is a weight filtration of $\dMgm(X)$: 
\[
M_{\le 1}(D)^*(2)[3] \in \DeffgM_{\BQ, w \le -1} \quad \text{and} \quad
M_{\le 1}(D) \in \DeffgM_{\BQ, w \ge 0} \; . 
\]
(d)~The isomorphism classes of the weight filtration from (c)
and of the Chow motive $M^{!*}(X^*)$ correspond under
the bijection from Theorem~\ref{1R} ($\BQ$-linear version).
\end{Thm}

\begin{Proof}
Parts~(a) and (b) follow from Proposition~\ref{2C}.
Part~(c) follows from stability under passage to direct factors~\ref{1I}~(1),
and the fact that
\[
M(D)^*(2)[3] \longto \dMgm(X) \longto M(D)
\longto M(D)^*(2)[4]
\]
is a weight filtration (Theorem~\ref{1O}).
Finally, (d) is a direct consequence of Construction~\ref{1Q},
and of the shape of the diagram~$\imm (\id - p)$.
\end{Proof}

\begin{Rem} \label{2E}
Let us discuss Construction~\ref{1Q} in the 
present geometrical setting.
The weight filtration of $\dMgm(X)$ is 
\[
M_{\le 1}(D)^*(2)[3] \stackrel{c_-}{\longto} \dMgm(X) 
\stackrel{c_+}{\longto} M_{\le 1}(D)
\stackrel{\delta}{\longto} M_{\le 1}(D)^*(2)[4] \; ;
\]
according to Theorem~\ref{2D}~(a), it is independent of the choice
of $\Xp$ (hence of $D$).
It fits into the diagram of exact triangles
\[
\vcenter{\xymatrix@R-10pt{
        0 &
        M_{\le 1}(D)^*(2)[4] \ar[l] &
        M_{\le 1}(D)^*(2)[4] \ar@{=}[l]  &
        0 \ar[l] \\
        \Mcgm(X) \ar[u] &
         &
        M_{\le 1}(D) \ar[u]_{\delta} &
        \Mcgm(X)[-1] \ar[u] \ar[l]_{c_+(v_+[-1])} \\
        \Mcgm(X) \ar@{=}[u]&  
        M(X) \ar[l]_-u &
        \dMgm(X) \ar[u]_{c_+} \ar[l]_-{v_-} &
        \Mcgm(X)[-1] \ar@{=}[u] \ar[l]_-{v_+[-1]} \\
        0 \ar[u] &
        M_{\le 1}(D)^*(2)[3] \ar[u]^{v_-c_-} \ar[l] &
        M_{\le 1}(D)^*(2)[3] \ar[u]_{c_-} \ar@{=}[l] &
        0 \ar[l] \ar[u]   
\\}}
\]
In the general situation of Construction~\ref{1Q}, what we did next
was to \emph{choose} some object $M_0$ completing the diagram.
In the specific situation we are considering at present, 
there is only one choice, up to unique isomorphism,
which in addition is compatible with 
any of the diagrams of type $(A)$ associated to desingularizations
$\Xp$ of $X^*$. This choice is $M^{!*}(X^*)$. We thus obtain rigidification
of the intersection motive, while the condition from Complement~\ref{1W}
on the absence of weights $-1$ and $0$ in the boundary motive 
is clearly not satisfied --- unless $X^* = X$ is itself (proper and)
smooth (cmp.~Problem~\ref{1U}).  
\end{Rem} 

Let us finish this section by an example, which will allow us
to illustrate both Principle~\ref{1B} (on extensions)
and Principle~\ref{1Ra} (on functoriality).

\begin{Ex} \label{2F}
Our base field $k$ equals the field of rational numbers $\BQ$.
Fix a real quadratic number field $L$, and let $X$
be a \emph{Hilbert modular surface} associated to $L$
and some level $K$. We view $K$ as an open compact subgroup
of the group $G(\BA_f)$ of (finite) adelic points of the group scheme $G$
from \cite[Sect.~1.27]{R}. The subgroup $K \subset G(\BA_f)$
is assumed to be sufficiently small, a condition which ensures
that $X$ is smooth over $\BQ$. Denote by $X^*$
its \emph{Baily--Borel compactification};
it is normal and projective over $\BQ$. \\[0.1cm]
(a)~The morphism
\[
i_*: M_{\le 1}(D) \longto M^{!*}(X^*)
\]
occurring in the weight filtration
\[
M_{\le 1}(D) \stackrel{i_*}{\longto} M^{!*}(X^*) 
\longto \Mcgm(X)
\longto M_{\le 1}(D)[1]
\]
of $\Mcgm(X)$ can be used to construct a morphism of certain
sub-quotients of its source and target,
\[
M_1(D) \longto M^{!*}_2(X^*)
\]
\cite[Thm.~6.6]{W4}. This morphism can be interpreted as an
element of
\[
\Ext^1_{DM^{eff}_{gm}(\BQ)_{\BQ}} 
\bigl( M_1(D)[-1] , M^{!*}_2(X^*)[-2] \bigr)  \; ,
\]
i.e., a one-extension in the triangulated category $DM^{eff}_{gm}(\BQ)_{\BQ}$
(note that according to \cite[Prop.~6.5]{W4}, $M_1(D)[-1]$ is an Artin motive).
Following \cite[Ex.~7.4]{W4},
it can be related to the \emph{Kummer--Chern--Eisenstein extensions}
considered in \cite{Cs}. In particular \cite[Ex.~7.4~(6)]{W4},
the extension is non-trivial. \\[0.1cm]
(b)~The intersection motive
$M^{!*}(X^*)$ carries a natural action of the Hecke algebra 
$R(K,G(\BA_f))$ associated to $K \subset G(\BA_f)$.
More precisely, let $x \in G(\BA_f)$. The Hilbert surface
$X$ is the target
of two finite \'etale morphisms $g_1,g_2: Y \to X$, where
$Y$ denotes the Hilbert surface of level $K' := K \cap x^{-1}Kx$.
In standard notation from the theory of Shimura varieties,
the morphism $g_1$ corresponds to $[\ \cdot 1]$, and the morphism $g_2$
to $[\ \cdot x^{-1}]$. Both morphisms can be extended to finite morphisms
between the Baily--Borel compactifications
\[
g_i : Y^* \to X^* \; ,
\]
satisfying the formulae $g_i^{-1}(X) = Y$, $i = 1,2$.
According to Theorem~\ref{2D}~(b), the diagram
\[
\vcenter{\xymatrix@R-10pt{
        0 &
        M_{\le 1}(D)^*(2)[4] \ar[l] &
        M_{\le 1}(D)^*(2)[4] \ar@{=}[l]  &
        0 \ar[l] \\
        \Mcgm(X) \ar[u] &
        M^{!*}(X^*) \ar[u]^{i^*} \ar[l] &
        M_{\le 1}(D) \ar[u] \ar[l]_{i_*} &
        \Mcgm(X)[-1] \ar[u] \ar[l] \\
        \Mcgm(X) \ar@{=}[u]&  
        M(X) \ar[u] \ar[l] &
        \dMgm(X) \ar[u] \ar[l] &
        \Mcgm(X)[-1] \ar@{=}[u] \ar[l] \\
        0 \ar[u] &
        M_{\le 1}(D)^*(2)[3] \ar[u] \ar[l] &
        M_{\le 1}(D)^*(2)[3] \ar[u] \ar@{=}[l] &
        0 \ar[l] \ar[u]   
\\}}
\]
is therefore stable under the composition $g_{2,*} g_1^*$.
By definition, this composition equals the action
of the class $KxK \in R(K,G(\BA_f))$. In particular,
the Hecke algebra acts on the whole of the above diagram.
It is useful to note that its effect on the third row, 
i.e., on the boundary triangle,
is the one described in Example~\ref{2bL}~(d), for $X= S$,
$\varphi = \id_Y$, and $e = \id_X$.
\end{Ex}

%%% Local Variables:
%%% mode: latex
%%% TeX-master: "head"
%%% End:

\bigskip

%\include{Sec3}

%%%%%%%%%%%%%%%%%%%%%%%%%%%%%%%%%%%%%%%%%%%%%%%%%%%%%%%%%%%%%%%%%%%%%%%
%
%  Section 3
%
%%%%%%%%%%%%%%%%%%%%%%%%%%%%%%%%%%%%%%%%%%%%%%%%%%%%%%%%%%%%%%%%%%%%%%%

\section{The interior motive of a product of universal elliptic curves}
\label{3}

%%%%%%%%%%%%%%%%%%%%%%%%%%%%%%%%

%%%%%%%%%%%%%%%%%%%%%%%%%%%%%%%%

In this section, our base field $k$ equals the field of rational numbers $\BQ$.
Fix integers $n \ge 3$ and $r \ge 0$, and let $S \in Sm/\BQ$
denote the modular curve parametri\-zing elliptic curves
with level $n$ structure.
Write $X \to S$ for the universal elliptic curve, and
\[
X^r := X \times_S \times \ldots \times_S X
\] 
for the $r$-fold fibre product of $X$ over $S$.  
Recall the decomposition 
\[
h(X/S) = \bigoplus_{i=0}^2 h_i(X/S)
\]
of the relative motive of $X$ from Theorem~\ref{2cA}.

\begin{Def} \label{3A}
Define ${}^r \CV \in \CHSM_\BQ$ as 
\[
{}^r \CV := 
\Sym^r h_1(X/S) \; .
\]
\end{Def}

The tensor product is in $\CHSM_\BQ$, and the symmetric powers are formed
with the usual convention concerning the (twist of) the natural
action of the symmetric group on a power of $X$ over $S$ (see e.g.\
\cite[Sect.~1.1.2]{Sch}). Thus, ${}^r \CV$ is a direct factor
of $h(X^r/S) \in \CHSM_\BQ$. 
That is, it is associated
to an idempotent
\[
e \in \ch^r (X^r \times_S X^r) \otimes_\BZ \BQ \; .
\]
From Theorem~\ref{2bB}~(a), we get a natural action of $e$
on the boundary triangle of $X^r$. In particular, we have the following
result.

\begin{Prop} \label{3B}
The triangle
\[
(\ast)_{X^r}^e \quad\quad
\dMgm(X^r)^e \longto M(X^r)^e \longto \Mcgm(X^r)^e \longto \dMgm(X^r)^e[1] 
\]
in $\DeffqgM$ is exact.
\end{Prop}

This triangle equals the image $( \dMgm , M , \Mcgm )(h(X^r/S)^e)$ 
of the smooth relative Chow motive ${}^r \CV = h(X^r/S)^e$ 
under the functor $( \dMgm , M , \Mcgm )_S$
from Theorem~\ref{2bB}.
It is not very difficult to check that the idempotent $e$
coincides with the one used in \cite[Sect.~1]{Sch} and \cite[Sect.~3 and 4]{W3}.

\begin{Prop}[{\cite[Ex.~4.16~(d)]{W3}}] \label{3C}
The direct factor $\dMgm(X^r)^e$ of the boundary motive of $X^r$
is without weights $-1$ and $0$ whenever $r \ge 1$.
\end{Prop}

By Complement~\ref{1W}, the $e$-part of the interior motive of $X^r$
can be constructed, and is unique up to unique isomorphism.
It is shown in \cite[Thm.~3.3~(b) and Cor.~3.4~(b)]{W3}
that it is canonically isomorphic to the Chow motive ${}^r_n \CW$ 
constructed in \cite{Sch} out of a compactification
of $X^r$. In that article, 
the action of the Hecke algebra on that compactification,
hence on ${}^r_n \CW$ is then used
to construct the Grothendieck motive $M(f)$ 
for elliptic normalized newforms $f$ of 
level $n \ge 3$ and weight $w = r + 2 \ge 3$.
Let us finish this section by giving an alternative
description of the action of the Hecke algebra
on ${}^r_n \CW$, which avoids compactifications.

\begin{Ex} \label{3D}
Assume that $r \ge 1$, and consider the diagram
\[
\vcenter{\xymatrix@R-10pt{
        0 &
        \dMgm(X^r)^e_{\le -2}[1] \ar[l] &
        \dMgm(X^r)^e_{\le -2}[1] \ar@{=}[l]  &
        0 \ar[l] \\
        \Mcgm(X^r)^e \ar[u] &
        {}^r_n \CW \ar[u] \ar[l]  &
        \dMgm(X^r)^e_{\ge 1} \ar[u] \ar[l] &
        \Mcgm(X^r)^e[-1] \ar[u] \ar[l] \\
        \Mcgm(X^r)^e \ar@{=}[u]&  
        M(X^r)^e \ar[u] \ar[l] &
        \dMgm(X^r)^e \ar[u] \ar[l] &
        \Mcgm(X^r)^e[-1] \ar@{=}[u] \ar[l] \\
        0 \ar[u] &
        \dMgm(X^r)^e_{\le -2} \ar[u] \ar[l] &
        \dMgm(X^r)^e_{\le -2} \ar[u] \ar@{=}[l] &
        0 \ar[l] \ar[u]   
\\}}
\]
from Complement~\ref{1W} associated to the weight filtration
of $\dMgm(X^r)^e$ avoiding weights $-1$ and $0$. We shall show that
this diagram, and hence ${}^r_n \CW$ in particular,
carries a natural action of the Hecke algebra 
$R(K_n,\GL_2(\BA_f))$ associated to the 
principal subgroup $K_n \subset GL_2(\BA_f)$ of level $n$.
Let $x \in \GL_2(\BA_f)$. The curve
$S$ is the target
of two finite \'etale morphisms $g_1,g_2: U \to S$, where
$U$ denotes the modular curve of level $K' := K_n \cap x^{-1}K_nx$.
In standard notation from the theory of Shimura varieties,
the morphism $g_1$ corresponds to $[\ \cdot 1]$, and the morphism $g_2$
to $[\ \cdot x^{-1}]$. 
Denote by $X_1, X_2$ the base changes of the universal
elliptic curve $X$ to $U$ \emph{via} $g_1$ and $g_2$,
respectively. To the data $K_n$ and $x$, the following are
canonically associated: a third elliptic curve $Y$ over $U$,
and isogenies $f_1: Y \to X_1$ and $f_2: Y \to X_2$.  
Now note that
$\varphi := \Gamma_{f_2^r} \circ {}^t \Gamma_{f_1^r}$
defines a morphism of smooth relative Chow motives over $U$,
\[
\varphi: h(X_1^r / U) = g_1^* ( h(X^r/S) ) \longto 
g_2^* ( h(X^r/S) ) = h(X_2^r / U) \; .
\] 
Since both $f_1$ and $f_2$ are isogenies, this morphism
is compatible 
with the external products of the idempotents $p_{X_i,1}$
projecting onto the $h_1$
(Theorem~\ref{2cA}~(b) and (c)). The morphism $\varphi$ is also 
compatible with the action of the symmetric group;
hence it is compatible with the cycle classes
$g_i^*(e_r)  
\in \ch^r (X^r_i \times_U X^r_i) \otimes_\BZ \BQ$. This means that
we have the relation
\[
\varphi \circ g_1^*(e_r) = g_2^*(e_r) \circ \varphi
\]
of morphisms of smooth relative Chow motives over $U$.
We are thus in the situation of Example~\ref{2bL}~(d),
and may therefore define the endomorphism $\varphi(g_1,g_2)$
of the boundary triangle $(\ast)_{X^r}^e$, i.e., of the third row
of the above diagram. The weight filtration 
of $\dMgm(X^r)^e$ being functorial, $\varphi(g_1,g_2)$
induces an endomorphism of the third column. 
Finally, thanks to Complement~\ref{1W},
the endomorphism extends uniquely to $M_0$. Altogether,
$\varphi(g_1,g_2)$ extends to the whole of the
above diagram of exact triangles. By definition,
this is the action of the class $K_n x K_n$ we aimed at.
\end{Ex}

% Local Variables:
%%% mode: latex
%%% TeX-master: "head"
%%% End:

\bigskip

%%%%%%%%%%%%%%%%%%%%%%%%%%%%%%%%%%%%%%%%%%%%%%%%%%%%%%%%%%%%%%%%%%%%%%%
%
%  Bibliography
%
%%%%%%%%%%%%%%%%%%%%%%%%%%%%%%%%%%%%%%%%%%%%%%%%%%%%%%%%%%%%%%%%%%%%%%%


\begin{thebibliography}{99}

\bibitem[A]{A}
Y.~Andr\'e, {\it Une introduction aux motifs}, 
Panoramas et Synth\`eses~{\bf 17}, Soc.\ Math.\ France (2004).

\bibitem[BBD]{BBD}
A.A.~Beilinson, J.~Bernstein, P.~Deligne, 
{\it Faisceaux pervers}, in: 
B.~Teissier, J.L.~Verdier (eds.), 
{\it Analyse et topologie sur les espaces singuliers (I)}, 
Ast\'erisque~{\bf 100}, Soc.\ Math.\ France (1982).

\bibitem[Bo1]{Bo1}
M.V.~Bondarko,
{\it Differential graded motives: weight complex, weight filtrations and 
spectral sequences for realizations; Voevodsky versus Hanamura},
J.\ Inst.\ Math.\ Jussieu~{\bf 8} (2009), 39--97.

\bibitem[Bo2]{Bo}
M.V.~Bondarko,
{\it Weight structures vs.\ $t$-structures; weight filtrations, 
spectral sequences, and complexes (for motives and in general)},
J.~$K$-Theory~{\bf 6} (2010), 387--504.

\bibitem[Cas]{Cs}
A.~Caspar, {\it Realisations of Kummer--Chern--Eisenstein-motives},
Manuscripta Math.~{\bf 122} (2007), 23--57.

\bibitem[CatMi]{CM}
M.A.A.~de Cataldo, L.~Migliorini, 
{\it The Chow motive of semismall re\-solutions},
Math.\ Res.\ Lett.~{\bf 11} (2004), 151--170.

\bibitem[CiD\'e]{CDeg}
D.-C.~Cisinski, F.~D\'eglise,
{\it Local and stable homological algebra
in Grothendieck abelian categories},
Homology, Homotopy and Applications~{\bf 11} (2009), 219--260.

\bibitem[D\'e1]{Deg1}
F.~D\'eglise,
{\it Finite correspondences and transfers over a regular base},
in J.~Nagel, C.~Peters (eds.), 
{\it Algebraic cycles and motives. Vol.~1. Proceedings of the EAGER Conference,
held in Leiden, August~30--September~4, 2004}, 
London Math.\ Soc.\ Lecture Note Ser.~{\bf 343}, 
Cambridge Univ.\ Press (2007), 138--205.

\bibitem[D\'e2]{Deg2}
F.~D\'eglise,
{\it Around the Gysin triangle II},
Doc.\ Math.~{\bf 13} (2008), 613--675.

\bibitem[DeMu]{DM}
C.~Deninger, J.~Murre,
{\it Motivic decomposition of abelian schemes and the Fourier transform},
J.\ reine angew.\ Math~{\bf 422} (1991), 201--219.

\bibitem[FV]{FV}
E.M.~Friedlander, V.~Voevodsky,
{\it Bivariant cycle cohomology},
Chapter~4 of \cite{VSF}.

\bibitem[H]{H}
G.~Harder,
{\it Eisensteinkohomologie und die Konstruktion gemischter Motive},
Lect.\ Notes Math.~{\bf 1562},
Springer-Verlag (1993).

\bibitem[L]{L3}
M.~Levine,
{\it Smooth Motives},
in R.~de Jeu, J.D.~Lewis (eds.),
{\it Motives and Algebraic Cycles. A Celebration in Honour of Spencer J.~Bloch},
Fields Institute Communications~{\bf 56}, 
American Math.\ Soc.\ (2009), 175--231.

\bibitem[R]{R}
M.~Rapoport,
{\it Compactification de l'espace de modules de Hilbert--Blumenthal},
Compositio Math.~{\bf 36} (1978), 255--335.

\bibitem[Sa]{Sa}
M.~Saito, {\it Mixed Hodge modules}, 
Publ.\ Res.\ Inst.\ Math.\ Sci.~{\bf 26} (1990), 221--333.

\bibitem[Sch]{Sch}
A.J.~Scholl, 
{\it Classical Motives}, 
in U.~Jannsen, S.~Kleiman, J.-P.~Serre (eds.), 
{\it Motives.
Proceedings of the AMS-IMS-SIAM Joint Summer Research Conference, 
held at the University of Washington, Seattle, July 20--August 2, 1991}, 
Proc.\ of Symp.\ in Pure Math.~{\bf 55}, Part 1, AMS (1994), 163--187.

\bibitem[St]{St}
J.~Steenbrink,
{\it A Summary of Mixed Hodge Theory},
in U.~Jannsen, S.~Kleiman, J.-P.~Serre (eds.), 
{\it Motives.
Proceedings of the AMS-IMS-SIAM Joint Summer Research Conference, 
held at the University of Washington, Seattle, July 20--August 2, 1991}, 
Proc.\ of Symp.\ in Pure Math.~{\bf 55}, Part 1, AMS (1994), 31--41.

\bibitem[V1]{V}
V.~Voevodsky,
{\it Triangulated categories of motives over a field},
Chapter~5 of \cite{VSF}.

\bibitem[V2]{V3}
V.~Voevodsky,
{\it Motivic cohomology groups are isomorphic to higher Chow groups 
in any characteristic},
Int.\ Math.\ Res.\ Notices~{\bf 2002} (2002), 351--355.

\bibitem[VSF]{VSF}
V.~Voe\-vodsky, A.~Suslin, E.M.~Friedlander,
{\it Cycles, Transfers, and Motivic Homology Theories},
Ann.\ of Math.\ Studies~{\bf 143}, Princeton Univ.\ Press (2000).

\bibitem[W1]{W1}
J.~Wildeshaus,
{\it The boundary motive: definition and basic properties},
Compositio Math.~{\bf 142} (2006), 631--656.

\bibitem[W2]{W2}
J.~Wildeshaus,
{\it On the boundary motive of a Shimura variety},
Compositio Math.~{\bf 143} (2007), 959--985.

\bibitem[W3]{W3}
J.~Wildeshaus,
{\it Chow motives without projectivity},
Compositio Math.~{\bf 145} (2009), 1196--1226.

\bibitem[W4]{W5}
J.~Wildeshaus,
{\it On the interior motive of certain Shimura varieties: 
the case of Hilbert--Blumenthal varieties},
preprint, June~2009, version dated March~18, 2011, 31~pages, submitted, 
available on arXiv.org under 
{\tt http://arxiv.org/abs/0906.4239} 

\bibitem[W5]{W4}
J.~Wildeshaus,
{\it Pure motives, mixed motives and extensions of motives
associated to singular surfaces}, 
39 pages, to appear in J.-B.~Bost and J.-M.~Fontaine (eds.), 
{\it Autour des motifs. Ecole d'\'et\'e franco-asiatique de g\'eom\'etrie 
alg\'ebrique et de th\'eorie des nombres. Vol.~II},
Panoramas et Synth\`eses, Soc.\ Math.\ France (2011), available on arXiv.org under 
{\tt http://xxx.lanl.gov/abs/0706.4447} 

\end{thebibliography}
\end{document}